\documentclass[leqno,11pt]{article}
\usepackage{latexsym}
\usepackage{amssymb}
\usepackage{amsfonts}
\usepackage{amsmath}
\usepackage{todonotes}
\usepackage{color}
\usepackage{url}
\usepackage{hyperref}

\definecolor{bleu}{rgb}{0.00,0.4,0.90}

\makeatletter
\long\def\unmarkedfootnote#1{{\long\def\@makefntext##1{##1}\footnotetext{#1}}}
\makeatother

\newtheorem{definition}{Definition}[section]
\newtheorem{lemma}[definition]{Lemma}
\newtheorem{theorem}[definition]{Theorem}
\newtheorem{proposition}[definition]{Proposition}

\newtheorem{remark}[definition]{Remark}

\def\Rvtex{\mathbb R}
\def\Nvtex{\mathbb N}
\def\rnvtex{{{\Rvtex}^n}}
\def\vavtex{\rightarrow}
\def\tovtex{\rightarrow}
\def\R{\mathbb R}

\def\rn{{{\R}^n}}

\def\r2{{{\R}^2}}

\def\va{\rightarrow}

\newcommand{\nequiv}{\not\equiv}

\def\bdm#1\edm{\begin{displaymath}#1\end{displaymath}}
\def\be#1\ee{\begin{equation}#1\end{equation}}
\def\barr#1\earr{\begin{align}#1\end{align}}

\newcommand{\medint}{-\kern  -,395cm\int}
\newcommand{\medintinrigo}{-\kern  -,315cm\int}
\newcommand{\medelle}{-\kern  -,235cm L}
\newcommand{\medellenrigo}{-\kern  -,180cm L}
\newcommand{\qed}{\thinspace\null\nobreak\hfill
\hbox{\vbox{\kern-.2pt\hrule height.2pt
depth.2pt\kern-.2pt\kern-.2pt \hbox to1.8mm {\kern-.2pt\vrule
width.4pt \kern-.2pt\raise1.8mm\vbox to.2pt{} \lower0pt\vtop
to.2pt{}\hfil\kern-.2pt \vrule
width.4pt\kern-.2pt}\kern-.2pt\kern-.2pt \hrule height.2pt depth.2pt
\kern-.2pt}}\par\medbreak}

\begin{document}

\title{An eigenvalue problem for the anisotropic $\Phi $-Laplacian}

\author{A. Alberico \thanks{Istituto per le
Applicazioni del Calcolo ``M. Picone'' (IAC),  Consiglio
Nazionale delle Ricerche (CNR), Via P. Castellino 111, 80131 Napoli,
Italy. E--mail:a.alberico@iac.cnr.it} -- G. di
Blasio\thanks{Dipartimento di Matematica e Fisica,
Universit\`{a} degli Studi della Campania \textquotedblleft L. Vanvitelli\textquotedblright , Viale Lincoln, 5 - 81100
Caserta, Italy. E--mail: giuseppina.diblasio@unicampania.it} -- F.
Feo\thanks{Dipartimento di Ingegneria, Universit\`{a} degli Studi di
Napoli \textquotedblleft Pathenope\textquotedblright, Centro
Direzionale Isola C4 80143 Napoli, Italy. E--mail:
filomena.feo@uniparthenope.it}}

\date{}

\maketitle

\begin{abstract}
We study an eigenvalue problem involving a fully anisotropic elliptic differential
operator in arbitrary Orlicz-Sobolev spaces. The relevant equations are
associated with constrained minimization problems for integral functionals
depending on the gradient of competing functions through general anisotropic
$N$-functions. In particular, the latter need neither be radial, nor have
a polynomial growth, and are not even assumed to satisfy the so called
$\Delta _{2}$-condition. The resulting analysis requires the development
of some new aspects of the theory of anisotropic Orlicz-Sobolev spaces.
\end{abstract}

\bigskip


\numberwithin{equation}{section}

\centerline{\textsc{How to cite this paper}}

\medskip

\noindent
This paper has been accepted to \emph{Journal of Differential Equations} and the final publication is available at


\centerline{\href{url}{https://doi.org/10.1016/j.jde.2020.03.049.}}

\medskip
\noindent
Should you wish to cite this paper, the authors would like to cordially ask you to cite it appropriately.

\bigskip

\section{\textbf{Introduction}}

In the present paper, we deal with the existence of solutions to  fully
anisotropic eigenvalue problems having the form
%
\begin{equation}
\label{EL}
\left \{
\begin{array}[c]{l@{\quad }l}
-\text{\textup{div}} \; (\Phi _{\xi }(\nabla u))= \lambda \, b(|u|) \,{\mathrm{sign}}
\, u & \qquad \text{\textup{in\ }} \Omega
\\
u=0 & \qquad \text{on}\;\partial \Omega \,,
\end{array}
\right .
\end{equation}
where $\Omega $ is an open bounded subset in $\rn $, with $n\geq 2$,
$\lambda $ is a positive real parameter,
$\Phi : \rn\va [ 0, \infty )$ is an $N$-function (see Section~\ref{sec_2.1}) belonging
to $\mathcal{C}^{1}(\rnvtex )$, and $b: [0, \infty )\vavtex [0, \infty )$ is an
increasing, left-continuous function such that $b(t)=0$ if and only if
$t=0$ and $\lim _{t\vavtex \infty } b(t)=+\infty $. Here, $\Phi _{\xi }$ denotes
the gradient of $\Phi $. Let us emphasize that $\Phi (\xi )$ neither necessarily
depends on $\xi $ through its length $|\xi |$, nor necessarily has a power
type behavior.
\par
Formally, problem {\eqref{EL}} represents the Euler-Lagrange equation associated
with the following constrained minimization problem
%
\begin{equation}
\label{min}
\inf \left \{  \int _{\Omega }\Phi (\nabla u)\, dx : u\in W^{1}_{0} L_{B,
\Phi } (\Omega ), \int _{\Omega }B(u)\, dx =r\right \}  \,,
\end{equation}
where $r$ is any positive real constant,
$W^{1}_{0}L_{B,\Phi } (\Omega )$ is the anisotropic Orlicz-Sobolev space
built upon $\Phi $ and $B$, where $B$ is the $1$-dimensional  $N$-function defined as
$B(t)= \int _{0}^{|t|} b(\tau )\, d\tau $. We point out that neither
$\Phi $ nor $B$ are require
to fulfill the
\emph{$\Delta _{2}$-condition}. Due to this fact, differentiability of the
functionals appearing in{\eqref{min}} is not guaranteed. Hence, problem {\eqref{EL}} cannot be derived via standard methods like constrained minimization
or critical point techniques.
\par
The function $B$ will be subject to a sharp growth condition that follows
from the anisotropic Sobolev inequality for
$W^{1}_{0}L_{B,\Phi } (\Omega )$ proved in \cite{C1}. For a comprehensive treatment of this matter, we refer
the reader to Section~\ref{sec2.3} and Section~\ref{sec3}.
\par
Our aim is to show that for any $r>0$ there exist $\lambda _{r}>0$ and
$u_{r}\in {W}_{0}^{1} L_{B,\Phi }(\Omega ){\cap L^{\infty } (\Omega )}$
such that $\int _{\Omega }B(u_{r})\,dx =r$ and $u_{r}$ solves problem {\eqref{EL}} with $\lambda =\lambda _{r}$.
\\
Classical results in this line of investigations deal with the eigenvalue
problem for $p$-Laplacian
%
\begin{equation}
\label{p_laplacian}
\left \{
\begin{array}[c]{l@{\quad }l}
-\text{\textup{div}} \; (|\nabla u|^{p-2}\nabla u)= \lambda |u|^{q-2}\,u &
\qquad \text{\textup{in\ }} \Omega
\\
u=0 & \qquad \text{on}\;\partial \Omega \,,
\end{array}
\right .
\end{equation}
with $1<p<N$ and $1<q<p^{*}$, where $p^{*}$ stands for the Sobolev conjugate
of $p$. Problem {\eqref{p_laplacian}} is the Euler-Lagrange equation associated
with the minimization problem {\eqref{min}} corresponding to the choice
$\Phi (\xi )= \frac{1}{p} |\xi |^{p}$. Several results are available in the
literature on existence and properties of eigenvalues and corresponding
eigenfunctions to problem {\eqref{p_laplacian}} (see, \emph{e.g.},
\cite{D,FNSS,L1,L2,SZ}).
\\
Isotropic eigenvalue problems and associated constrained minimization problems
in the spirit of {\eqref{EL}} and {\eqref{min}}, respectively, with
$\Phi (\xi )= \Phi (|\xi |)$ and $B(t)= \Phi (|\xi |)$, have been investigated
in \cite{MT}. Our contribution extends the results of \cite{MT}, not only
in allowing for completely fully anisotropic differential operators, but
also in admitting more general growths on the right-hand side
$b(|u|) \,{\mathrm{sign}}\, u$. In particular, the generality of the problems
under consideration calls for the use and further development of the unconventional
functional framework of anisotropic Orlicz and Orlicz-Sobolev spaces which
are not necessarily reflexive (see, \emph{e.g.},
\cite{BC,CGZ-G, C1,C2,C3,C4, Sc,Sk1,Sk2,Tr}).
\\
Let us mention that elliptic equations and variational problems, whose
growth is governed by an $n$-dimensional $N$-function $\Phi $, have been
studied under different perspectives in
\cite{A,AdBF1,AdBF2,AC,ACCZ-G,BC,Ch,C1,C2,C3,C4, GWWZ}.
\\
The paper is organized as follows. Section~\ref{sec2} contains a background, as
well as some new results, on anisotropic Orlicz and Orlicz-Sobolev spaces.
The statements of our main results and some special instances are given
in Section~\ref{sec3}. The proofs of main results are presented in Section~\ref{sec4}.


\section{Functional setting}
\label{sec2}

\subsection{Young functions}
\label{sec_2.1}

Let $n\geq 1$. Let
$\Phi :\mathbb{R}^{n}\rightarrow \left [ 0,+\infty \right ]$ be an
$n$-\textit{dimensional Young function}, namely an even, convex function
such that $\Phi (0)=0$ and, for every $t>0$, the set
$\{\xi \in \rnvtex : \Phi (\xi )<t\}$ is bounded and contains an open neighborhood
of $0$. An $n$-dimensional Young function $\Phi $ is called an
\emph{$n$-dimensional $N$-function} if it is a finite valued function, vanishes
only at $0$ and the following additional conditions are in force
%
\begin{equation}
\underset{\left \vert \xi \right \vert \rightarrow +\infty }{\lim }
\frac{\Phi \left (
\xi \right ) }{\left \vert \xi \right \vert }=+\infty
\label{lim_inf}%
\end{equation}
and %
%
\begin{equation}
\underset{\left \vert \xi \right \vert \rightarrow 0}{\lim }
\frac{\Phi \left (
\xi \right ) }{\left \vert \xi \right \vert }=0\,.
\label{lim_0}
\end{equation}

For $n=1$, any \emph{$1$-dimensional $N$-function}
$A:\mathbb{R}\vavtex [0,+\infty )$ takes the form
%
\begin{equation}
\label{B}
A(t)= \int _{0}^{|t|} a(\tau ) \; d\tau \qquad \text{for}\;\; t\in \Rvtex ,
\end{equation}
where $a: [0, \infty ) \vavtex [0, \infty )$ is an increasing, right-continuous
function, which is positive for $\tau >0$ and satisfies conditions
$a(0)=0$ and $\lim _{\tau \rightarrow +\infty } a(\tau )= +\infty $.

If $\Phi $ is an $n$-dimensional Young function, then
%
\begin{equation}
\label{lin}
\,\Phi (h\xi )\leq |h|\Phi (\xi ) \qquad \text{for}\;\; |h|\leq 1 \;
\text{and}\; \xi \in \rnvtex \,.
\end{equation}
The Young inequality tells us that
%
\begin{equation}
\xi \cdot \xi ^{\prime }\leq \Phi \left ( \xi \right ) +\Phi _{
\bullet }\left ( \xi ^{\prime }\right ) \ \ \text{ for }\ \xi ,\xi ^{
\prime }\in \mathbb{R}%
^{n},
\label{Young}%
\end{equation}
where $\Phi _{\bullet }$ is the \textit{Young conjugate} of $\Phi $ given
by%
%
\begin{equation}
\Phi _{\bullet }\left ( \xi ^{\prime }\right ) =\sup \left \{  \xi
\cdot \xi ^{\prime }-\Phi \left ( \xi \right ) :\xi \in \mathbb{R}^{n}
\right \}  \ \ \text{ for }\ %
\xi ^{\prime }\in \mathbb{R}^{n}.
\label{Young_func}%
\end{equation}
Here, \textquotedblleft\,$\cdot$\,\textquotedblright\ stands for scalar
product in ${{\mathbb{R}}^{n}}$. We observe that if
$\Phi $ is finite-valued and assumption {\eqref{lim_inf}} holds, then the
function $\Phi _{\bullet }$ is an $n$-dimensional Young function and finite-valued
(see \cite[Corollary 6.3]{BC}). Note also that the Young conjugation is involutive,
\emph{i.e.} $\Phi _{\bullet \bullet }= \Phi $. Moreover,
$\Phi _{\bullet }$ is an $N$-function, provided that $\Phi $ is.
\\
An $n$-dimensional Young function $\Phi $ is said to satisfy the
\emph{$\Delta _{2}$-condition} near infinity, briefly
$\Phi \in \Delta _{2}$ near infinity, if it is finite-valued
and there exist constants $C>2$ and $K\geq 0$ such that
$\Phi (2\xi )\leq C\, \Phi (\xi )$ for $|\xi |>K$.
\\
Let us consider a case when the $n$-dimensional $N$-function $\Phi $ is
given by
%
\begin{equation}
\label{ex_Ai}
\Phi (\xi )= \sum _{i=1}^{n} A_{i}(\xi _{i})\qquad \text{for
$\xi \in \rn $,}
\end{equation}
where $A_{i}$, for $i=1,\ldots ,n$, are $1$-dimensional $N$-functions.
A standard choice in {\eqref{ex_Ai}} is $A_{i}(t)=|t|^{p_{i}}$ for some powers
$1<p_{i}<+\infty $, for $i = 1, \dots , n$. One can easily verify that
in {\eqref{ex_Ai}} every function $A_{i}\in \Delta _{2}$ near infinity if
and only if $\Phi (\xi )$ does. An example of a function which does not
satisfy the $\Delta _{2}$-condition is given by
\begin{equation*}
\Phi (\xi )=\sum _{i=1}^{n} (e^{|\xi _{i}|^{\alpha _{i}}}-1)\qquad
\text{for $\xi \in \rnvtex $}
\end{equation*}
with $\alpha _{i}>1$, for any $i=1,\ldots ,n$.
\\
The following proposition is a special case of
\cite[Theorem 5.1]{Sk1}.
%
\begin{proposition}[Equality cases in the Young inequality]
\label{equalityYoung}
Let $\Phi $ be a differentiable
$n$-dimensional Young function. Then, for any $\xi _{0}\in \rnvtex $
\begin{equation*}
\xi _{0} \cdot \eta =\Phi \left ( \xi _{0} \right ) + \Phi _{\bullet }
\left ( \eta \right )
\end{equation*}
if and only if $\eta = \Phi _{\xi }(\xi _{0})$.
\end{proposition}

Thanks to {Proposition~\ref{equalityYoung}}, in
\cite[Proposition 6.7]{BC} the authors proved the following lemma when
$\Phi \in \mathcal{C}^{1}(\rnvtex )$, but their proof runs also under the weaker assumption
that $\Phi $ is differentiable.

\begin{lemma}%
\label{lemmaBC}%
 \textup{\textbf{\cite{BC}}} Let $\Phi $ be a differentiable $n$-dimensional
Young function. Assume that {\eqref{lim_inf}} holds. Then
%
\begin{equation}
\label{dis_gradiente}
\Phi _{\bullet }(\Phi _{\xi }(\xi ))\leq \Phi _{\xi }(\xi ) \cdot \xi
\leq \Phi (2\xi )\qquad \text{for $\xi \in \rnvtex $}\, .
\end{equation}
\end{lemma}

Finally, we show a technical lemma which will be very useful in the sequel.
\\
We say that two $n$-dimensional $N$-functions $\Phi $ and $\Psi $ are equivalent
if there exist positive constants $k_{1}$ and $k_{2}$, depending only on
$n$, such that
\begin{equation*}
\label{equiv}
\Phi (k_{1}\xi ) \leq \Psi (\xi )\leq \Phi (k_{2}\xi )\quad \text{for}
\;\; \xi \in \rnvtex .
\end{equation*}
We emphasize that $\Phi $ and $\Psi $ are equivalent if and only if
$\Phi _{\bullet }$ and $\Psi _{\bullet }$ are.

\begin{lemma}%
\label{young_covex}
Given any $n$-dimensional $N$-function $\Phi $, there exists another
$n$-dimensional $N$-function $\Psi $ which is strictly convex and equivalent
to $\Phi $. As a consequence, $\Psi _{\bullet }$ is differentiable.
\end{lemma}

\medskip
\noindent \emph{Proof.}
\rm
Theorem $26.3$ in \cite{Ro} states that the strict
convexity of an $N$-function guarantees the differentiability of its
conjugate. Thus, it is enough to prove the existence of a strictly convex
$N$-function equivalent to $\Phi $. Let $\Phi _{-}: \rnvtex \tovtex [0,
\infty )$ be the radial function defined as
\begin{equation*}
\Phi _{-}(\xi )= \sup \Big \{\Theta (\xi ): \;\;\Theta :\rnvtex \tovtex [0,
\infty ) \;\; {\text{$N$-function, radial and
$\Theta (\xi )\leq \Phi (\xi )$}}\Big \} \qquad \text{for}\;\;\xi \in
\rnvtex \,.
\end{equation*}
By construction, $\Phi _{-}$ is an $N$-function. Indeed,
$\Phi _{-}$ is a convex function since it is a supremum of convex functions,
and one can easily check that conditions {\eqref{lim_inf}} and {\eqref{lim_0}} are verified.
\\
Fixed $c>0$ and let $g:[0, \infty ) \tovtex [0,\infty )$ be a strictly increasing
function such that $0<g(s)\leq c$ for $s\geq 0$. Then,


%
\begin{equation*}
G(t)=\int _{0}^{t} g(s)\; ds
\end{equation*}
is a strictly convex, increasing function and $0<G(t)\leq c\,t $ for any
$t>0$. Set
\begin{equation*}
\Upsilon (\xi )=G(\Phi _{-}(\xi ))\,.
\end{equation*}
Since $\Phi _{-}$ is radial and $G$ is strictly convex, it follows that
also $\Upsilon $ is strictly convex.
\\
Then, $\Phi + \Upsilon $ is an $N$-function, strictly convex and equivalent
to $\Phi $ because
\begin{equation*}
\Phi (\xi ) \leq \Phi (\xi )+ \Upsilon (\xi )\leq \Phi (\xi ) + c
\Phi _{-}(\xi )\leq (1+c)\Phi (\xi ) \leq \Phi \big ((1+c)\xi \big )\,,
\end{equation*}
where the last inequality is due to {\eqref{lin}}.

\qed

\subsection{Anisotropic Orlicz spaces}
\label{sec2.2}

In this section we present Orlicz spaces built upon both a $1$-dimensional
Young function (see, \emph{e.g.}, \cite{Adams}) and $n$-dimensional Young
functions (see, \emph{e.g.}, \cite{BC,Sk1,Sk2,Sc}). For the convenience
of the reader we give a briefly background.
\\
Let $\Omega $ be a bounded measurable subset in $\rnvtex $, with
$n\geq 2$. The \emph{Orlicz space} $L_{A} (\Omega )$, associated with a
$1$-dimensional Young function $A$, is the set of all measurable functions
$g: \Omega \tovtex \Rvtex $ such that the Luxemburg norm
\begin{equation*}
\|g\|_{L_{A}(\Omega )}= \inf \left \{  k >0 : \int _{\Omega }A \left (
\frac{g(x)}{k} \right ) dx \leq 1 \right \}
\end{equation*}
is finite. The functional $\|\, \cdot \, \|_{L_{A}(\Omega )}$ is a norm
on $L_{A}(\Omega )$, which makes the latter a Banach space.
\\
Given two finite-valued $1$-dimensional Young functions $A$ and $D$, we
say that $A\prec \prec D$, namely $A$
\emph{increases essentially more slowly} than $D$ near infinity, if
\begin{equation*}
\label{esse_slow}
\lim _{t\vavtex +\infty } \frac{A(\gamma t)}{D(t)}=0 \qquad \text{for every}
\;\;\gamma >0\,.
\end{equation*}
Note that if $A\prec \prec D$, then
\begin{equation*}
L_{D}(\Omega ) \hookrightarrow L_{A}(\Omega )\, ,
\end{equation*}
where the arrow ``\,$\hookrightarrow$\,'' stands for continuous embedding.

Let $\Phi $ be an $n$-dimensional Young function. The anisotropic
\emph{Orlicz class} $\mathcal{L}_{\Phi }(\Omega ; \rnvtex )$ is defined as
\begin{equation*}
\label{orlicz_class}
\mathcal{L}_{\Phi }(\Omega ;\mathbb{R}^{n})=\left \{  U:\Omega
\rightarrow \mathbb{R}^{n} \;\text{measurable s.t.}\; {\displaystyle
\int _{\Omega }} \Phi \left ( U\right ) \; dx<+\infty \right \}  \,.
\end{equation*}
Note that $\mathcal{L}_{\Phi }(\Omega ; \rnvtex )$ is a convex set of functions
and it needs not be a linear space in general, unless $\Phi $ satisfies
the $\Delta _{2}$-condition near infinity. The \emph{Orlicz space}
$L_{\Phi }(\Omega ; \rnvtex )$ is the linear hull of
$\mathcal{L}_{\Phi }(\Omega ; \mathbb{R}^{n})$ and it is a Banach space
with respect to the following Luxemburg norm
%
\begin{equation}
\label{Lux_norm}
\left \Vert U\right \Vert _{\Phi }=\inf \left \{  k>0:%
{\displaystyle \int _{\Omega }} \Phi \left ( \frac{U}{k}\right )\; dx
\leq 1\right \}  \,.
\end{equation}
We emphasize that
$L_{\Phi }(\Omega ; \rnvtex )\subset L^{1}(\Omega )$ for any $n$-dimensional
Young function $\Phi $. We stress that if two $n$-dimensional Young functions
$\Phi $ and $\Psi $ are equivalent, then $\|\cdot \|_{\Phi }$ and
$\|\cdot \|_{\Psi }$ are equivalent and then $L_{\Phi }$ and $L_{\Psi }$ are
the same space.
\\
Let us denote by $E_{\Phi }(\Omega ; \rnvtex )$ the closure in
$L_{\Phi }(\Omega ; \rnvtex )$ of the bounded measurable functions with compact
support in $\overline{\Omega }$. In general
%
\begin{equation}
\label{600}
E_{\Phi }(\Omega ; \rnvtex )\subset \mathcal{L}_{\Phi }(\Omega ;\mathbb{R}^{n})
\subset L_{\Phi }%
(\Omega ;\mathbb{R}^{n})\,.
\end{equation}
Both inclusions hold as equalities in {\eqref{600}} if and only if
$\Phi $ satisfies the $\Delta _{2}$-condition near infinity (see
\cite[Corollary 5.1]{Sc}).

From now on, let $\Phi $ be an $n$-dimensional $N$-function. The following
generalized H\"{o}lder inequality holds
%
\begin{equation}
\label{holder}
\int _{\Omega } U(x)\cdot V(x)\;dx\leq 2\, \| U\|_{\Phi }\| V\|_{
\Phi _{\bullet }}
\end{equation}
for every $U\in L_{\Phi }(\Omega ; \mathbb{R}^{n})$ and
$V\in L_{\Phi _{\bullet }}(\Omega ; \mathbb{R}^{n})$ (see
\cite[Theorem
4.1]{Sk2}). Fixed $V\in L_{\Phi _{\bullet }}(\Omega ; \rnvtex )$, the integral
in {\eqref{holder}} defines a linear and continuous functional on
$L_{\Phi }(\Omega ; \mathbb{R}^{n})$. The space
$L_{\Phi }(\Omega ;\mathbb{R}^{n})$ can be also endowed with the following
Orlicz norm
%
\begin{equation}
\label{orlicz_norm}
\| U\|_{(\Phi )}=
\underset{\int _{\Omega }\Phi _{\bullet }(V)\leq 1}{\sup }\bigg | {
\displaystyle \int _{\Omega }} U(x)\cdot V(x)\;dx\bigg |\,.
\end{equation}
Thanks to {Lemma~\ref{young_covex}}, we can assume that
$\Phi $ is differentiable (up to an equivalent $N$-function), and then
Luxemburg norm {\eqref{Lux_norm}} and Orlicz norm {\eqref{orlicz_norm}} are
equivalent, \textit{i.e.}
$\|U\|_{\Phi }\leq \|U\|_{(\Phi )}\leq 2\|U\|_{\Phi }$ (see
\cite[Theorem 4.5]{Sk2}).
\\
Combining the Orlicz norm and the Luxemburg norm together it is possible
to get this sharp form of generalized H\"{o}lder inequality
\begin{equation*}
\label{Holder_2}
{\displaystyle \int _{\Omega }} U(x)\cdot V(x)\;dx\leq \left \Vert U
\right \Vert _{(\Phi )}\left \Vert V\right \Vert _{\Phi _{\bullet }}
\end{equation*}
for every $U\in L_{\Phi }(\Omega ; \mathbb{R}^{n})$ and
$V\in L_{\Phi _{\bullet }}(\Omega ; \mathbb{R}^{n})$.
\smallskip

If $A$ is a $1$-dimensional $N$-function, it is well-known that the dual
space of $E_{A}(\Omega )$ is isomorphic and homeomorphic to
$L_{A_{\bullet }}(\Omega )$ (see \cite[Theorem 8.18]{Adams}). The analogue
result holds for the anisotropic spaces.

\begin{proposition}%
\label{Prop_dual}
Let $\Phi $ be an $n$-dimensional $N$-function. The dual space of
$E_{\Phi }(\Omega ;\mathbb{R}^{n})$, denoted by
$(E_{\Phi }(\Omega ;\rnvtex ))^{\prime }$, is isomorphic and homeomorphic
to $L_{\Phi _{\bullet }}(\Omega ;\mathbb{R}^{n})$ and the duality pairing
is given by
\begin{equation*}
<V,U>= {\displaystyle \int _{\Omega }} V(x)\cdot U(x)dx
\end{equation*}
for $V\in L_{\Phi _{\bullet }}(\Omega ;\mathbb{R}^{n})$ and
$U\in E_{\Phi }%
(\Omega ;\mathbb{R}^{n})$.
\end{proposition}

%
\begin{remark}%
\label{remark2.5}
Note that if $\Phi \in \Delta _{2}$, then
$\left ( L_{\Phi }(\Omega ;\mathbb{R}%
^{n})\right ) ^{\prime }=L_{\Phi _{\bullet }}(\Omega ;\mathbb{R}^{n})$.
\end{remark}

\medskip

\textbf{Proof of Proposition \ref{Prop_dual}.}
\rm
We proceed by steps. First
we show that any element
$V\in L_{\Phi _{\bullet }}(\Omega ;\mathbb{R}^{n})$ determines a bounded
linear functional defined as
%
\begin{equation}
\label{pair}
<\textbf{\emph{l}}_{V},U>=\int _{\Omega } U(x)\cdot V(x) dx
\end{equation}
for every $U\in E_{\Phi }(\Omega ;\mathbb{R}^{n})$. Then, it remains to show
that every bounded linear functional on $E_{\Phi }(\Omega ;\mathbb{R}^{n})$
can be written uniquely in the form $\textbf{\emph{l}}_{V}$ for some $V \in
L_{\Phi _{\bullet }}(\Omega ;\mathbb{R}^{n})$. In order to do this, we prove
that any bounded linear functional $\textbf{\emph{l}}$ on $E_{\Phi }(\Omega
;\mathbb{R}^{n})$ has the form {\eqref{pair}} when we restrict ourselves to the
set of vector-valued simple functions, \textit{i.e.} vector-valued functions
such that each component is a simple function (functions that assume a finite
number of values). The density of this set in $E_{\Phi }$ allows us to
conclude the proof.

\textbf{Step 1.} \textit{$\textbf{l}_{V}$ restricted on
$E_{\Phi }(\Omega ;\mathbb{R}^{n})$ belongs to
$\left (E_{\Phi }(\Omega ;\mathbb{R}^{n})\right )'$}.

It follows by {\eqref{holder}}.

\textbf{Step 2.} \textit{The set of vector-valued simple functions is
dense in $E_{\Phi }$.}

Let us consider $U\in L^{\infty }(\Omega ;\mathbb{R}^{n})$. By standard
measure theory, there exists a sequence $\{U_{h}\}_{h}$ of vector-valued
simple functions such that

i) $U_{h}\rightarrow U$ \textit{a.e.} in $\Omega $,

ii)
$\|U_{h}\|_{L^{\infty }(\Omega ;\mathbb{R}^{n})}\leq C_{0}\|U\|_{L^{
\infty }(\Omega ;\mathbb{R}^{n})}$ $\forall h \in \Nvtex $ and for a given
positive constant $C_{0}$.

We claim that $\|U_{h}-U\|_{\Phi }\rightarrow 0$. Indeed, by ii) we have
\begin{equation*}
\|U_{h}-U\|_{L^{\infty }(\Omega ;\mathbb{R}^{n})}\leq (C_{0}+1)\|U\|_{L^{
\infty }(\Omega ;\mathbb{R}^{n})}:=C_{1}.
\end{equation*}
Given $\epsilon >0$, set
\begin{equation*}
C_{2}(\epsilon )= \underset{|\xi |\leq C_{1}}{\sup } \Phi \left (
\frac{\xi }{\epsilon }\right )
\end{equation*}
and
\begin{equation*}
S_{h}=\Big \{x\in \Omega:\Phi \left (\frac{U_{h}(x)-U(x)}{\epsilon }\right )>
\frac{1}{2|\Omega |}\Big \}.
\end{equation*}
By i) and continuity of $\Phi $, the measure $|S_{h}|\rightarrow 0$ as
$h\rightarrow \infty $. Then, we can choose $h$ sufficiently large such
that $|S_{h}|<\frac{1}{2C_{2}(\epsilon )}$. It follows that
\begin{align*}
\int _{\Omega } \Phi \left (\frac{U_{h}-U}{\epsilon }\right )\,dx &=
\int _{S_{h}} \Phi \left (\frac{U_{h}-U}{\epsilon }\right )\,dx+\int _{
\Omega \setminus S_{h}} \Phi \left (\frac{U_{h}-U}{\epsilon }\right )
\,dx
\\
&\leq C_{2}(\epsilon ) \frac{1}{2C_{2}(\epsilon )}+|\Omega |
\frac{1}{2|\Omega |} = 1,
\nonumber
\end{align*}
that implies $\|U_{h}-U\|_{\Phi }\leq \epsilon $. By definition of
$E_{\Phi }(\Omega ;\mathbb{R}^{n})$, the result holds for
$U \in E_{\Phi }(\Omega ;\mathbb{R}^{n})$.

\textbf{Step 3.} \textit{Representation formula for $\textbf{l}\in (E_{\Phi
}(\Omega ;\mathbb{R}^{n}))'$ on vector-valued simple functions.}

Let $\textbf{\emph{l}}$ be any continuous linear functional on
$E_{\Phi }(\Omega ;\mathbb{R}^{n})$ and let $G\subset \Omega $ be a measurable
set. Let us consider $\chi _{G}(x)$, the characteristic function of the
set $G$, and let us set
\begin{equation*}
{\overrightarrow{\chi }}_{G}(x)=(
\underbrace{\chi _{G}(x),\cdots ,\chi _{G}(x)}_{n}) \text{ and }
\overrightarrow{\chi } _{G,h}(x)=(\underbrace{0,\ldots , 0}_{h-1},
\chi _{G}^{h}(x),\underbrace{0,\ldots ,0}_{n-h})\,.
\end{equation*}
Then
${\overrightarrow{\chi }}_{G}(x) =\sum _{h=1}^{n}\chi _{G }(x)e_{h}=
\sum _{h=1}^{n}\overrightarrow{\chi } _{G,h}(x)$, where
$\{e_{1},\ldots , e_{n}\}$ is the standard base of $\mathbb{R}^{n}$. For
every $h\in \{1,\ldots ,n\}$, by Step 2 we have that
$\overrightarrow{\chi } _{G,h}\in E_{\Phi }(\Omega ;\mathbb{R}^{n})$. Moreover
$\mu _{h}(G)=<\textbf{\emph{l}},\overrightarrow{\chi } _{G,h}(x)>$ is an
absolutely continuous measure. Indeed,
%
\begin{equation}
\label{ass_cont_1}
|<\textbf{\emph{l}},\overrightarrow{\chi } _{G,h}(x)>|\leq \|
\textbf{\emph{l}}\| \quad \|\overrightarrow{\chi } _{G,h}\|_{\Phi }
\end{equation}
with
\begin{equation*}
\|\overrightarrow{\chi } _{G,h}\|_{\Phi }=\inf \Big \{k:\int _{\Omega }
\Phi \left (\frac{\overrightarrow{\chi }
_{G,h}(x)}{k}\right )\, dx\leq 1\Big \}=\inf \left \{  k:|G|\Phi \left (
\frac{e_{h}}{k}\right )\leq 1\right \}  .
\end{equation*}
Let us define the $1$-dimensional Young function
$A_{h}: \mathbb{R}\rightarrow [0,\infty )$ as
$A_{h}(t)=\Phi (t e_{h})$ for every $t\in \mathbb{R}$. Then,
%
\begin{equation}
\label{ass_cont_2}
\|\overrightarrow{\chi } _{G,h}\|_{\Phi }=
\frac{1}{A^{-1}_{h}\left (\frac{1}{|G|}\right )}.
\end{equation}
The absolute continuity of measure $\mu _{h}$ follows by combining {\eqref{ass_cont_1}} and {\eqref{ass_cont_2}}. By virtue of the Radon-Nikodym's
Theorem, there exists a real valued function $V_{h}$ belonging to
$L^{1}(\Omega )$ such that
%
\begin{equation}
\label{302}
<\textbf{\emph{l}},\overrightarrow{\chi } _{G,h}>=\int _{\Omega}V_{h}(x)
\chi _{G}(x) \,dx \quad \forall h=1,\dots ,n.
\end{equation}
By {\eqref{302}},
%
\begin{equation}
\label{33}
<\textbf{\emph{l}},\overrightarrow{\chi } _{G}>=\sum _{h=1}^{n}<
\textbf{\emph{l}},\overrightarrow{\chi } _{G,h}>=\int _{\Omega}V(x)\cdot
\overrightarrow{\chi } _{G}(x) \,dx\,,
\end{equation}
where $V(x)=(V_{1}(x),\dots ,V_{n}(x))$. Moreover, if $U$ is a vector-valued
simple function defined as
\begin{equation*}
\label{semplice}
U(x)=\sum _{j=1}^{n} \alpha _{j} \overrightarrow{\chi }_{G_{j}}(x),
\end{equation*}
where $\alpha _{j}\in \mathbb{R}$ and $G_{j}$ are disjoint measurable subsets
of $\Omega $, by the linearity of $\textbf{\emph{l}}$ and {\eqref{33}}, we
get
%
\begin{equation}
\label{l}
<\textbf{\emph{l}},U>=\sum _{j=1}^{n} \alpha _{j} <\textbf{\emph{l}},
\overrightarrow{\chi }_{G_{j}}(x)>= \sum _{j=1}^{n}\alpha _{j} \int _{
\Omega } V(x)\cdot \overrightarrow{\chi }_{G_{j}}(x)\, dx= \int _{
\Omega } U \cdot V \,dx= <\textbf{\emph{l}}_{V},U>,
\end{equation}
where $\textbf{\emph{l}}_{V}$ is defined by {\eqref{pair}}.

\textbf{Step 4.} \textit{Function $V$ that appears in {\eqref{l}} belongs
to $ L_{\Phi _{\bullet }}(\Omega ;\mathbb{R}^{n})$}.

Let $U\in E_{\Phi }(\Omega ;\mathbb{R}^{n})$. By Step 2, we know that there
exists a sequence $U_{h}$ of simple functions such that
$U_{h} \rightarrow U$ in $L_{\Phi }(\Omega ;\mathbb{R}^{n})$. This means
that $U_{h} \rightarrow U$ almost everywhere and also the sequence
$|U_{h}\cdot V|\rightarrow |U\cdot V|$ almost everywhere. Moreover, fixed
some positive constant $K$, one can choose $h$ sufficiently large such
that
\begin{equation*}
\|U_{h}\|_{\Phi }\leq \|U\|_{\Phi }+\|U_{h}-U\|_{\Phi }\leq \|U\|_{
\Phi }+K.
\end{equation*}
Now, if $U_{h} \cdot V\geq 0$, on applying Fatou's Lemma, we get
%
\begin{equation}
\label{val_ass_1}
\int _{\Omega } U \cdot V\, dx\leq \liminf \int _{\Omega } U_{h}
\cdot V\, dx\leq \liminf \bigg |\int _{\Omega } U_{h} \cdot V \,dx
\bigg | \leq \liminf \|\textbf{\emph{l}}\| \|U_{h}\|_{\Phi }\leq \|
\textbf{\emph{l}}\|(\|U\|_{\Phi }+K).
\end{equation}
On the other hand, if $U_{h} \cdot V\leq 0$, on applying Fatou's Lemma
again, we get
%
\begin{align}
\label{val_ass_2}
\int _{\Omega } -U \cdot V\, dx&\leq \liminf \int _{\Omega } -U_{h}
\cdot V\, dx\leq \liminf \bigg |\int _{\Omega } -U_{h} \cdot V \,dx
\bigg | \leq \liminf \|\textbf{\emph{l}}\| \|U_{h}\|_{\Phi }
\nonumber
\\
&\leq \|
\textbf{\emph{l}}\|(\|U\|_{\Phi }+K).
\end{align}
By {\eqref{val_ass_1}} and {\eqref{val_ass_2}}, we deduce
%
\begin{equation}
\label{55}
\bigg |\int _{\Omega } U \cdot V\, dx\bigg |< +\infty
\end{equation}
for any $U\in E_{\Phi }(\Omega ;\mathbb{R}^{n})$. This means that if we
choose $U(x)=(U_{1}(x),\dots ,U_{n}(x))$ such that
\begin{equation*}
U_{i}(x)=\frac{\partial \Phi _{\bullet }}{\partial e_{i}}(x) \quad
\text{for}\;\; i=1,\dots ,n
\end{equation*}
then, by {Proposition~\ref{equalityYoung}}, we get
\begin{equation*}
\int _{\Omega } \Phi _{\bullet }(V) \, dx\leq \int _{\Omega } \Phi _{
\bullet }(V)\, dx+\int _{\Omega } \Phi (U)\, dx=\int _{\Omega } U
\cdot V \,dx\,.
\end{equation*}
We stress that the extra assumption on differentiability of
$\Phi _{\bullet }$ required in {Proposition~\ref{equalityYoung}} can be dropped
thanks to {Lemma~\ref{young_covex}}. Finally, since {\eqref{55}}, it follows
that $V\in L_{\Phi _{\bullet }}(\Omega ;\mathbb{R}^{n})$ and
$\textbf{\emph{l}}_{V}$ is linear bounded functional on
$E_{\Phi }(\Omega ,\mathbb{R}^{n})$.

\textbf{Step 5.} \textit{Identification between} $\textbf{\emph{l}}$
\textit{and} $\textbf{\emph{l}}_{V}$.

We note that both the functionals $\textbf{\emph{l}}_{V}$ defined as in {\eqref{pair}} and $\textbf{\emph{l}}$ assume the same values on the set of
vector-valued simple functions. Since the last set is dense in
$E_{\Phi }(\Omega ,\mathbb{R}^{n})$, they agree with
$E_{\Phi }(\Omega ,\mathbb{R}^{n})$ and the proof is complete.

\qed

For the convenience of the reader, let us recall a few definitions concerning
the convergence and boundedness. We say that a sequence
$\{U_{h}\}_{h}\subset L_{\Phi }(\Omega ;\mathbb{R}^{n})$
\emph{converges in mean} to $U\in L_{\Phi }(\Omega ;\mathbb{R}^{n})$ if
\begin{equation*}
\int _{\Omega }\Phi (U_{h} -U)\; dx \longrightarrow 0\,.
\end{equation*}
Note that the convergence in norm implies the convergence in mean, and
they are equivalent if and only if $\Phi $ satisfies the $\Delta _{2}$-condition.
Moreover, if a function $U\in L_{\Phi }(\Omega ;\mathbb{R}^{n})$ is
\emph{bounded in mean}, namely
\begin{equation*}
\int _{\Omega }\Phi (U)\; dx < C\,
\end{equation*}
for some constant $C>0$, then it is bounded in norm as well. The converse
is not true unless $\Phi $ satisfies the $\Delta _{2}$-condition.

\qed
\bigskip

Finally, we extend \cite[Lemma 1]{Go2} to vector-valued functions that
it will be useful in what follows.

\begin{lemma}%
\label{lemma_sup}
Let $\Phi $ be a differentiable $n$-dimensional $N$-function. For all
$V\in L_{\Phi }(\Omega ; \mathbb{R}^{n}) $
\begin{align*}
&\sup \left \{  \int _{\Omega }U \cdot V \; dx - \int _{\Omega }\Phi _{\bullet }(U)\; dx \;: \;U\in L_{\Phi _{\bullet }}(\Omega ; \mathbb{R}^{n})
\right \}  = \int _{\Omega }\Phi (V)\; dx
\\
&= \sup \left \{  \int _{\Omega }U \cdot V \; dx - \int _{\Omega }\Phi _{\bullet }(U)\; dx \;: \;U\in E_{\Phi _{\bullet }}(\Omega ; \mathbb{R}^{n})
\right \}
\nonumber
\end{align*}
\end{lemma}

\medskip
\noindent \emph{Proof.}
\rm
Since
$E_{\Phi _{\bullet }}(\Omega ; \mathbb{R}^{n})\subset L_{\Phi _{\bullet }}(
\Omega ; \mathbb{R}^{n})$, we have only to prove that
%
\begin{align}
\label{intermedia}
&\sup \left \{  \int _{\Omega }U \cdot V \; dx - \int _{\Omega }\Phi _{\bullet }(U)\; dx \;: \;U\in L_{\Phi _{\bullet }}(\Omega ; \mathbb{R}^{n})
\right \}  \leq \int _{\Omega }\Phi (V)\; dx
\\
&\leq \sup \left \{  \int _{\Omega }U \cdot V \; dx - \int _{\Omega }\Phi _{\bullet }(U)\; dx \;: \;U\in E_{\Phi _{\bullet }}(\Omega ; \mathbb{R}^{n})
\right \}  .
\nonumber
\end{align}
The left-hand side in {\eqref{intermedia}} follows by applying Young inequality.

Now let us prove the right-hand side in {\eqref{intermedia}}. Let
$V\in L_{\Phi }(\Omega ; \mathbb{R}^{n})$. We define
$U_{h}=\Phi _{\xi }(V_{h})$ by
\begin{equation*}
V_{h}=\left \{
\begin{array}[c]{l@{\quad }l}%
V(x) & \text{for $x$ \emph{s.t.} $|V(x)|\leq h$}
\\
\\
0 & \text{otherwise,}
\end{array}
\right .
\label{P_sym}%
\end{equation*}
where $|\cdot |$ denotes the modulus in $\mathbb{R}^{n}$. We claim that
$U_{h}\in L^{\infty }(\Omega ; \mathbb{R}^{n})\subset E_{\Phi _{\bullet }}(\Omega ; \mathbb{R}^{n}) $. Indeed, by {\eqref{dis_gradiente}} we get
%
\begin{equation}
\label{88}
0\leq \Phi _{\bullet }(U_{h})\leq \Phi (2V_{h}).
\end{equation}
By definition of $N$-function and {\eqref{88}}, it easily follows that if
$V\equiv 0$ then $U_{h}\equiv 0$ and the claim is obvious. Let
$V\nequiv 0$ and suppose by contradiction that
$\sup _{x\in \Omega }|U_{h}(x)|=+\infty $. Then, there exists a sequence
$\{x_{j}\}_{j}\subset \Omega $ such that
$|U_{h}(x_{j})|\vavtex +\infty $. Thus, by {\eqref{88}}, we get
\begin{equation*}
0\leq
\frac{\Phi _{\bullet }\left (U_{h}(x_{j})\right )}{\left |U_{h}(x_{j})\right |}
\leq
\frac{\Phi \left (2V_{h}(x_{j})\right )}{\left |U_{h}(x_{j})\right |\,.}
\end{equation*}
By {\eqref{lim_inf}}, the term in the center blows up, while the right-hand
side goes to zero, because $\Phi $ is continuous and $V_{h}(x)$ is bounded.
This proves that $U_{h}\in L^{\infty }(\Omega ; \mathbb{R}^{n})$.

By {Proposition~\ref{equalityYoung}}, we obtain
$\Phi (V_{h})= U_{h} \cdot V_{h}-\Phi _{\bullet }(U_{h})$. So, by integrating
on $\Omega $, it follows that
\begin{equation*}
\int _{\Omega }\Phi (V_{h})\;dx= \int _{\Omega }U_{h} \cdot V_{h} \;dx -
\int _{\Omega }\Phi _{\bullet }(U_{h}) \;dx \leq \int _{\Omega }U_{h} \cdot V
\;dx -\int _{\Omega }\Phi _{\bullet }(U_{h})\;dx
\end{equation*}
and then
%
\begin{equation}
\label{303}
\int _{\Omega }\Phi (V_{h})\;dx\leq \sup \left \{  \int _{\Omega }U \cdot V
\; dx - \int _{\Omega }\Phi _{\bullet }(U)\; dx \;: \;U\in E_{\Phi _{\bullet }}(\Omega ; \mathbb{R}^{n}) \right \}  .
\end{equation}
By Fatou's Lemma, the left-hand side in {\eqref{303}} converges to
$\int _{\Omega }\Phi (V)\;dx<+\infty $ and {\eqref{intermedia}} follows.

\qed

\subsection{Anisotropic Orlicz-Sobolev spaces}
\label{sec2.3}

Let $\Phi $ be an $n$-dimensional $N$-function. Let us define the Banach
space $\mathcal{W}^{1}_{0} L_{\Phi }(\Omega )$ (see \cite{BC}) as
\begin{align}
\mathcal{W}^{1}_{0} L_{\Phi }(\Omega ) = \{ u:\Omega \vavtex \Rvtex :\; &\text{ the
continuation of u by} \;\;0 \;\; \text{outside} \;\; \Omega
\nonumber
\\
& \;\;\text{is weakly differentiable and}\; \; \nabla u\in L_{\Phi }(
\Omega ; \rnvtex )\}
\nonumber
\end{align}
equipped with the norm
\begin{equation*}
\label{norm_W_Andrea}
\left \|  u\right \|  _{\mathcal{W}^{1}_{0}L_{\Phi }(\Omega )}=\left \|
\nabla u\right \|  _{L_{\Phi }(\Omega ;\rnvtex )}\,.
\end{equation*}
We emphasize that the following anisotropic Sobolev type inequality holds
for any function in $\mathcal{W}^{1}_{0} L_{\Phi }(\Omega )$ (see
\cite{C1}). Assume that $\Phi $ fulfils
%
\begin{equation}
\label{cond_in_zero}
\int _{0} \left (\frac{\tau }{\Phi _{\circ }(\tau )}\right )^{
\frac{1}{n-1}} \; d\tau <\infty ,
\end{equation}
where $\Phi _{\circ }: [0, \infty ) \vavtex [0, \infty )$ is an $N$-function
satisfying
\begin{equation*}
|\{\xi \in \rnvtex : \Phi (\xi )\leq t\}|= |\{\xi \in \rnvtex : \Phi _{\circ }(|
\xi |)\leq t\}| \qquad \text{for} \;\; t\geq 0.
\end{equation*}
Note that the function $\xi \vavtex \Phi _{\circ }(|\xi |)$ agrees with the spherically
increasing symmetral of $\Phi $.
\\
We denote by $\Phi _{n}:[0, \infty )\vavtex [0, \infty ]$ the optimal Sobolev
conjugate of $\Phi $ defined as
%
\begin{equation}
\label{sobolev_conj}
\Phi _{n}(t)= \Phi _{\circ }(H^{-1}(t)) \qquad \text{for}\;\; t\geq 0,
\end{equation}
where $H:[0, \infty )\vavtex [0, \infty )$ is given by
\begin{equation*}
\label{H}
H(t)= \left ( \int _{0}^{t} \left (\frac{\tau }{\Phi _{\circ }(\tau )}
\right )^{\frac{1}{n-1}} \; d\tau \right )^{\frac{n-1}{n}}\qquad \text{for}
\;\; t\geq 0,
\end{equation*}
provided that the integral is convergent. Here, $H^{-1}$ denotes the generalized
left-continuous inverse of $H$.

If
%
\begin{equation}
\label{int_infinito}
\int ^{\infty }\left (\frac{\tau }{\Phi _{\circ }(\tau )}\right )^{
\frac{1}{n-1}} \; d\tau =\infty \,,
\end{equation}
then there exists a constant $C_{1}=C_{1}(n)$ such that
%
\begin{equation}
\label{Opt_Sobolev_emb}
\|u\|_{L_{\Phi _{n}}(\Omega )}\leq C_{1}\| u\|_{\mathcal{W}_{0}^{1}L_{\Phi }(\Omega )}
\end{equation}
for every $u\in \mathcal{W}^{1}_{0} L_{\Phi }(\Omega )$ (see
\cite[Theorem 1 and Remark 1]{C1}).
\goodbreak

If
%
\begin{equation}
\label{int_finito}
\int ^{\infty }\left (\frac{\tau }{\Phi _{\circ }(\tau )}\right )^{
\frac{1}{n-1}} \; d\tau <\infty \,,
\end{equation}
then there exists a constant $C_{2} = C_{2}(\Phi , n, |\Omega |)$ such
that
%
\begin{equation}
\label{Opt_Sobolev_emb_infty}
\|u\|_{L^{\infty }(\Omega )}\leq C_{2}\| u\|_{\mathcal{W}_{0}^{1}L_{\Phi }(\Omega)}
\end{equation}
for every $u\in \mathcal{W}^{1}_{0} L_{\Phi }(\Omega )$ (see
\textit{e.g.} \cite[Theorem 1.2]{C4}).

\medskip

We define the anisotropic Orlicz-Sobolev space
$W^{1}L_{B,\Phi }(\Omega )$ as
\begin{equation*}
\label{def_W^1L_Phi}
W^{1}L_{B,\Phi }(\Omega ) =\left \{  u\in L_{B}(\Omega ): u\;\; \text{is weakly
differentiable in $\Omega $ and}\;\;\nabla u\in L_{\Phi }(\Omega ;
\mathbb{R}^{n})\right \}  \,.
\end{equation*}
The space $W^{1}E_{B,\Phi }(\Omega )$ is defined accordingly by replacing
$L_{B}(\Omega )$ and $L_{\Phi }(\Omega ;\rnvtex )$ by $E_{B}(\Omega )$ and
$E_{\Phi }(\Omega ;\rnvtex )$, respectively. Both
$W^{1}L_{B,\Phi }(\Omega )$ and $W^{1}E_{B,\Phi }(\Omega )$ can be identified
to subspaces of the product
$L_{B}(\Omega )\times L_{\Phi }(\Omega ; \rnvtex )$. The spaces
$W^{1}L_{B,\Phi }(\Omega )$ and $W^{1}E_{B,\Phi }(\Omega )$ equipped with
the norm
%
\begin{equation}
\label{norm_W_bis}
\left \Vert u\right \Vert _{W^{1}L_{B,\Phi }(\Omega )}=\left \Vert u
\right \Vert _{L^{B}(\Omega )}+\left \Vert \nabla u\right \Vert _{L_{
\Phi }(\Omega ;\rnvtex )}
\end{equation}
are Banach spaces (see \cite[Theorem 3.2]{Adams}).

Let us denote by $W^{1}_{0} L_{B,\Phi }(\Omega )$ the
$\sigma (L_{B}\times L_{\Phi }, E_{B_{\bullet }}\times E_{\Phi _{
\bullet }})$-closure of $\mathcal{D}(\Omega )$ in
$W^{1}L_{B,\Phi }(\Omega )$. Analogously,
$W^{1}_{0} E_{B,\Phi }(\Omega )$ stands for the closure of
$\mathcal{D}(\Omega )$ in $W^{1}L_{B,\Phi }(\Omega )$ with respect to
the norm {\eqref{norm_W_bis}}.

Let us emphasize that, given a function
$u\in W^{1}_{0} L_{B,\Phi }(\Omega )$, the function obtained by extending
$u$ outside $\Omega $ by zero belongs to $ W^{1}L_{B,\Phi }(\rnvtex )$. Thus,
%
\begin{equation}
\label{W_subset_W}
W^{1}_{0}L_{B,\Phi }(\Omega ) \subset \mathcal{W}^{1}_{0} L_{\Phi }(
\Omega )\,.
\end{equation}
Both spaces, $W_{0}^{1}L_{B,\Phi }(\Omega )$ and
$\mathcal{W}^{1}_{0} L_{\Phi }(\Omega )$, are reflexive if and only if
$\Phi \in \Delta _{2}$ near infinity.
\\
Embedding {\eqref{W_subset_W}} yields directly that Sobolev type inequalities {\eqref{Opt_Sobolev_emb}} and {\eqref{Opt_Sobolev_emb_infty}} hold for
$W^{1}_{0}L_{B,\Phi }(\Omega )$.
\\
Moreover, the following compact embedding holds.

\begin{proposition}%
\label{Prop_compact_embed}
Let $\Phi $ be an $N$-function. Assume that either {\eqref{int_infinito}} holds and $B\prec \prec \Phi _{n}$ or {\eqref{int_finito}} holds and $B$ is anything. Then
%
\begin{equation}
\label{compact_imbedding_1}
W^{1}_{0} L_{B,\Phi }(\Omega ) \hookrightarrow \hookrightarrow E_{B}(
\Omega )\,,
\end{equation}
where the  arrows ``$\, \hookrightarrow\hookrightarrow\,$'' stand for compact
embedding.
\end{proposition}

\medskip
\noindent \emph{Proof.}
\rm
Arguing as in the proof of \cite[Theorem
2.1]{CS}, we deduce that
\begin{equation*}
\label{500_bis}
\bigg \{ u: \int _{\Omega }\Phi _{n}\left (\frac{|u(x)|}{\lambda }\right )
\; dx <\infty \; \text{for every}\;\lambda >0\bigg \} \subset \text{closure
of }\; L^{\infty }(\Omega ) \; \text{in}\; L_{\Phi _{n}}(\Omega )\,.
\end{equation*}
Finally, observing that $\int _{\Omega }\Phi _{n}\left (\frac{|u(x)|}{\lambda
}\right )\; dx < \infty $ for every $\lambda >0$ whenever $u\in
W_{0}^{1}L_{B,\Phi }(\Omega )$ (see also \cite[Remark 7]{C0}), we conclude
that, if $B\prec \prec \Phi _{n}$,
%
\begin{equation}
\label{5000}
W_{0}^{1}L_{B,\Phi }(\Omega ) \subset E_{\Phi _{n}}(\Omega )\subset E_{B}(
\Omega )\,.
\end{equation}
\\
Let $\{u_{h}\}_{h}$ be a bounded sequence in
$W^{1}_{0} L_{B,\Phi }(\Omega ) $. Since the compact embedding (see
\cite{BC})
\begin{equation*}
\label{compact_imbedding_2}
W^{1}_{0} L_{B,\Phi }(\Omega ) \vavtex W^{1,1}_{0}(\Omega )
\hookrightarrow \hookrightarrow L^{1}(\Omega )\,,
\end{equation*}
it follows that (up a subsequence) $\{u_{h}\}_{h}$ converges in
$L^{1}(\Omega )$ and then in measure in $\Omega $. If
$B\prec \prec \Phi _{n}$, the convergence in measure
and the boundedness in $L_{\Phi _{n}}(\Omega )$ of $\{u_{h}\}_{h}$ (that
follows by {\eqref{Opt_Sobolev_emb}}) yield that $\{u_{h}\}_{h}$ converges in
$L_{B}(\Omega )$ (see \cite[Theorem 8.22]{Adams}). The embedding {\eqref{5000}} and the closure of $E_{B}(\Omega )$ conclude the proof.

If {\eqref{int_finito}} holds, we have that for any $N$-function
$B$,
%
\begin{equation}
\label{A100}
W_{0}^{1}L_{B,\Phi }(\Omega ) \subset L^{\infty }(\Omega)\subset E_{B_{1}}(
\Omega ) \subset E_{B}(\Omega )\,,
\end{equation}
where $B_{1}$ is an $N$-function such that $B\prec \prec B_{1}$. Thus,
one can use the same argument as in the previous case on replacing
$\Phi _{n}$ by $B_{1}$ in order to conclude the proof.

\qed

\subsection{Complementary systems}
\label{sec_compl_syst}

Let $X$ and $K$ be real Banach spaces in duality with respect to continuous
pairing $\langle \, \cdot , \cdot \,\rangle $, and let $X_{0}$ and
$K_{0}$ be subspaces of $X$ and $K$, respectively. Then,
$\left (X, X_{0}; K, K_{0}\right )$ represents a so-called
\emph{complementary system} if, by means of
$\langle \, \cdot , \cdot \,\rangle $, the dual of $X_{0}$ can be identified
to $K$ and that of $K_{0}$ to $K$.
\\
Given a complementary system $\left (X, X_{0}; K, K_{0}\right )$ and a
closed subspace $Y$ of $X$, it is possible to construct a new complementary
system imposing some restrictions on $Y$. More precisely, set
$Y_{0}= Y\cap X_{0}$, $Z= K/Y_{0}^{\perp }$ and
$Z_{0}=\{z+ Y_{0}^{\perp }: z\in K_{0}\}\subset Z$, where
$Y_{0}^{\perp }=\{ z\in K: \langle y, z \rangle =0\;\; \text{for every}
\; y\in Y_{0}\}$.
\\
The theory on complementary system has been investigated \emph{e.g.} in
\cite{Go1}, and, for the convenience of the reader, we recall Lemma 1.2
contained in it. The relevant lemma provides conditions so that
$\left (Y, Y_{0}; Z, Z_{0}\right )$ is a complementary system generated
by $Y$ in $\left (X, X_{0}; K, K_{0}\right )$.

\begin{lemma}%
\label{lemma_1.2_Gossez}
The pairing $\langle \, \cdot , \cdot \,\rangle $ between $X$ and
$K$ induces a pairing between $Y$ and $Z$ if and only if $Y_{0}$ is
$\sigma (X, K)$ dense in $Y$. In this case,
$\left (Y, Y_{0}; Z, Z_{0}\right )$ is a complementary system if $Y$ is
$\sigma (X, K_{0})$ closed, and conversely, when $K_{0}$ is complete,
$Y$ is $\sigma (X, K_{0})$ closed if
$\left (Y, Y_{0}; Z, Z_{0}\right )$ is a complementary system.
\end{lemma}
The topologies $\sigma (Y, Z)$ and $\sigma (Y, Z_{0})$ are the weak topologies
induced on $Y$ by $\sigma (X, K)$ and $\sigma (X, K_{0})$, respectively,
and $Z_{0}$ is the subspace of the dual space of $Y_{0}$ equals the set
of those linear functionals on $Y_{0}$ which are
$\sigma (X_{0}, K_{0})$ continuous.

Here, our aim is to prove that $Y=W^{1}_{0}L_{B,\Phi }(\Omega )$ generates
a new complementary system in
$(X,X_{0};K,K_{0})=\left (L_{B}\times L_{\Phi }, E_{B}\times E_{\Phi }; L_{B_{\bullet }}\times L_{\Phi _{\bullet }}, E_{B_{\bullet }}\times E_{\Phi _{\bullet }}\right )$.
\\
In order to do this, we assume that $\Omega $ enjoys the
\emph{segment property}, namely there exist a locally finite open covering
$\{\Omega _{j}\}_{j}$ of $\partial \Omega $ and corresponding vectors
$\{y_{j}\}_{j}$ such that $x+ty_{j}\in \Omega $ with
$x\in \overline{\Omega } \cap \Omega _{j}$ and $0<t<1$. This condition
is essential in {Lemma~\ref{lemma_E_dense}} below.

Let us verify that the conditions in {Lemma~\ref{lemma_1.2_Gossez}} are fulfilled.
First, $W^{1}_{0}L_{B,\Phi }(\Omega )$ is
$\sigma (L_{B}\times L_{\Phi }, E_{B_{\bullet }}\times E_{\Phi _{\bullet }})$
closed thanks to the very definition of
$W^{1}_{0}L_{B,\Phi }(\Omega )$. Moreover, we have to verify that
$W^{1}_{0}L_{B,\Phi }(\Omega )\cap (E_{B}\times E_{\Phi })$ agrees with
$W^{1}_{0}E_{B,\Phi }(\Omega )$ and it is
$\sigma (L_{B}\times L_{\Phi }, L_{B_{\bullet }}\times L_{\Phi _{\bullet }})$
dense in $W^{1}_{0}L_{B,\Phi }(\Omega )$.

\begin{lemma}%
\label{lemma_E_dense}
If $\Omega $ enjoys the segment property, then
\begin{itemize}
\item[\textup{\textbf{(a)}}] $W^{1}_{0}E_{B,\Phi }(\Omega )= W^{1}_{0}L_{B,
\Phi }(\Omega )\cap (E_{B}\times E_{\Phi })$,
\item[\textup{\textbf{(b)}}] $W^{1}_{0}E_{B,\Phi }(\Omega )$ is
$\sigma (L_{B}\times L_{\Phi }, L_{B_{\bullet }}\times L_{\Phi _{\bullet }})$
dense in $W^{1}_{0}L_{B,\Phi }(\Omega )$.
\end{itemize}
\end{lemma}

\medskip
\noindent \emph{Proof.}
\rm
\textbf{(a)} To prove the assert, one can be reduced
to prove that $\mathcal{D}(\Omega )$ is norm dense in $W^{1}_{0}L_{B,\Phi
}(\Omega )\cap (E_{B}\times E_{\Phi })$. By \cite[Theorem 1.3 and Corollary
1.10]{Go1}, it is enough to verify that $\mathcal{D}(\Omega )$ is $\sigma
(L_{B}\times L_{\Phi }, L_{B_{\bullet }}\times L_{\Phi _{\bullet }})$ dense
in $W^{1}_{0}L_{B,\Phi }(\Omega )$. It follows by  an appropriate
version of Lemmas 1.4 - 1.7 in \cite{Go1} applied to the $n$-dimensional
Young function $\Phi $. In fact, one can easily verify that those lemmas
hold for vectorial functions, as well.
\smallskip

\noindent
\textbf{(b)} Let us recall
that, by definition, $\mathcal{D}(\Omega )$ is dense in $W^{1}_{0}E_{B,\Phi
}(\Omega )$ with respect to the norm {\eqref{Lux_norm}} and is $\sigma
(L_{B}\times L_{\Phi }, E_{B_{\bullet }}\times E_{\Phi _{\bullet }})$ dense
in $W^{1}_{0}L_{B,\Phi }(\Omega )$. Our goal is to prove that
$W^{1}_{0}E_{B,\Phi }(\Omega )$ is $\sigma (L_{B}\times L_{\Phi },
L_{B_{\bullet }}\times L_{\Phi _{\bullet }})$ dense in $W^{1}_{0}L_{B,\Phi
}(\Omega )$, namely that for every $u\in W^{1}_{0}L_{B,\Phi }(\Omega )$ there
exists a sequence $\{u_{h}\}_{h}\subset W_{0}^{1}E_{B,\Phi }(\Omega )$ such
that
\begin{equation*}
\label{398}
u_{h} \vavtex u \qquad \text{in}\;\; \sigma (L_{B}\times L_{\Phi }, L_{B_{\bullet }}\times L_{\Phi _{\bullet }})\,,%
\end{equation*}
\emph{i.e.}
\begin{equation*}
\label{399}
\int _{\Omega }u_{h} \psi _{1} \; dx + \int _{\Omega }\nabla u_{h} \cdot
\psi _{2}\; dx \longrightarrow \int _{\Omega }u \psi _{1} \; dx + \int _{\Omega }\nabla u\cdot \psi _{2}\; dx \qquad \forall (\psi _{1}, \psi _{2})
\in L_{B_{\bullet }} \times L_{\Phi _{\bullet }}\,.
\end{equation*}

Let us suppose by contradiction that there exists a function
$\bar{u}$ in $W_{0}^{1}L_{B,\Phi }(\Omega )$ such that, for every sequence
$\{u_{h}\}_{h}\subset W_{0}^{1}E_{B,\Phi }(\Omega )$,
%
\begin{equation}
\label{400}
\lim _{h\vavtex \infty } u_{h} \neq \bar{u}\qquad \text{ in
$\sigma (L_{B}\times L_{\Phi }, L_{B_{\bullet }}\times L_{\Phi _{\bullet }})$}
\,.
\end{equation}
On the other hand, by the very definition of
$W_{0}^{1}E_{B,\Phi }(\Omega )$, for every
$u_{h}\in W_{0}^{1}E_{B,\Phi }(\Omega )$ there exists a sequence
$\{v_{h}^{k}\}_{k\in \Nvtex } \subset \mathcal{D} (\Omega )$ such that
$ v_{h}^{k} \vavtex u_{h}$ in norm, and then
%
\begin{equation}
\label{402}
\lim _{k\vavtex \infty }v_{h}^{k}= u_{h} \qquad \text{in}\;\; \sigma (L_{B}
\times L_{\Phi }, L_{B_{\bullet }}\times L_{\Phi _{\bullet }})\,.
\end{equation}
The statement follows by observing that {\eqref{402}} does
not agree with {\eqref{400}} and by recalling that $W^{1}_{0} L_{B,\Phi }(\Omega )$ is the
$\sigma (L_{B}\times L_{\Phi }, E_{B_{\bullet }}\times E_{\Phi _{
\bullet }})$-closure of $\mathcal{D}(\Omega )$ in
$W^{1}L_{B,\Phi }(\Omega )$.

\qed

\smallskip
{Lemma~\ref{lemma_1.2_Gossez}} and {Lemma~\ref{lemma_E_dense}} assure that
$\bigl(W^{1}_{0} L_{B,\Phi }(\Omega ), W^{1}_{0} E_{B,\Phi }(\Omega ); W^{-1}
L_{B_{\bullet },\Phi _{\bullet }}(\Omega ),\break
W^{-1}E_{B_{\bullet },\Phi _{\bullet }}(\Omega )\bigr)$ is the complementary
system generated by $W^{1}_{0} L_{B,\Phi }(\Omega )$ in
$(L_{B}\times L_{\Phi }$,
$ E_{B}\times E_{\Phi }; L_{B_{\bullet }}\times L_{\Phi _{\bullet }}, E_{B_{\bullet }}\times E_{\Phi _{\bullet }})$, where
\begin{equation*}
\label{dual_space_W}
W^{-1}L_{B_{\bullet },\Phi _{\bullet }}(\Omega )\,{=}\left \{  f\in
\mathcal{D}'(\Omega ): f\,{=}\,f_{0}-\sum _{i=1}^{n}
\frac{\partial }{\partial x_{i}} f_{i},\; \; f_{0}\,{\in}\, L_{B_{\bullet }}(
\Omega ), (f_{1},...,f_{n}) \,{\in}\, L_{\Phi _{\bullet }}(\Omega ;
\mathbb{R}^{n}) \right \}
\end{equation*}
and
\begin{align*}
W^{-1}E_{B_{\bullet },\Phi _{\bullet }}(\Omega )&=\left \{  f\in
\mathcal{D}'(\Omega ): f=f_{0}-\sum _{i=1}^{n}
\frac{\partial }{\partial x_{i}} f_{i},\; \; f_{0}\in E_{B_{\bullet }}(
\Omega ), (f_{1},...,f_{n}) \in E_{\Phi _{\bullet }}(\Omega ;\rnvtex )
\right \}  \,.
\end{align*}

\section{Main results}
\label{sec3}

Assume that $\Omega $ is an open bounded subset of $\mathbb{R}^{n}$, with
$n\geq 2$, satisfying the segment property, and
$\Phi \in \mathcal{C}^{1}(\rnvtex )$ is an $N$-function fulfilling {\eqref{cond_in_zero}}.

Our first main result concerns the existence of solutions to the following
Dirichlet problem
%
\begin{equation}
\label{EL_bis}
\left \{
\begin{array}[c]{l@{\quad }l}
-\text{\textup{div}} \; (\Phi _{\xi }(\nabla u))= \lambda \, \,{b(|u|)}
\text{ sign} \,u & \qquad \text{\textup{in\ }} \Omega
\\
u=0 & \qquad \text{on}\;\partial \Omega ,
\end{array}
\right .
\end{equation}
where $\Phi _{\xi }$ denotes the gradient of $\Phi $, $\lambda >0$ and for
$t>0$ function $b(t)$ is the derivative (see {\eqref{B}}) of a $1$-dimensional
$N$-function $B$ fulfilling some suitable assumptions.

%
\begin{definition}%
\label{def_weaksol}
A function $u\in \mathcal{W}_{0}^{1}L_{\Phi }(\Omega )$ is called a
\emph{weak solution} to problem {\eqref{EL_bis}} if
$\Phi _{\xi }(\nabla u)\in L_{\Phi _{\bullet }}(\Omega ; \rnvtex )$,
$b(|u|)\in L_{B_{\bullet }}(\Omega )$ and
%
\begin{equation}
\label{def_sol}
\int _{\Omega }\Phi _{\xi }(\nabla u) \cdot \nabla \varphi \; dx =
\lambda \int _{\Omega }\frac{b(|u|)}{|u|}u \varphi \; dx
\end{equation}
for any
$\varphi \in \mathcal{W}_{0}^{1}L_{\Phi }(\Omega )\cap L^{\infty }(
\Omega )$.
\end{definition}%

The existence result for solutions to {\eqref{EL_bis}} reads as follows.

%
\begin{theorem}%
\label{main_th}
Let $\Omega $ be a bounded Lipschitz domain in $\rnvtex $. Let
$\Phi \in \mathcal{C}^{1}(\rnvtex )$ be an $n$-dimensional $N$-function fulfilling {\eqref{cond_in_zero}}. Assume that $B$ is a $1$-dimensional $N$-function
such that $B\prec \prec \Phi _{n}$ if {\eqref{int_infinito}} holds or
$B$ is any if {\eqref{int_finito}} is in force. Then, for any $r>0$ there
exists a weak solution
$u_{r}\in \mathcal{W}_{0}^{1} L_{\Phi }(\Omega ){\cap L^{\infty } (
\Omega )}$ to problem {\eqref{EL_bis}}, with $\lambda =\lambda _{r}$, such
that $\int _{\Omega }B(u_{r})\,dx =r$.
\end{theorem}%
%

%
\begin{remark}%
\label{remark2}
A close inspection of the proof of Theorem~\ref{main_th} reveals that the
solution $u_{r}$ actually belongs to the space
${W}_{0}^{1} L_{B,\Phi }(\Omega )\subset \mathcal{W}_{0}^{1} L_{\Phi }(
\Omega )$ (see {\eqref{W_subset_W}}).
\end{remark}%
%

%
\begin{remark}%
\label{remark3.4}
Note that any bounded Lipschitz domain in $\rnvtex
$ satisfies the segment property, also.
\end{remark}%

In order to establish our main result we consider the following constrained
minimization problem
%
\begin{equation}
\label{min_pbl}
c_{r}= \inf \left \{  \int _{\Omega }\Phi (\nabla u)\;\; dx : u\in W^{1}_{0}
L_{B,\Phi } (\Omega ), \int _{\Omega }B(u)\;\;dx =r\right \}
\end{equation}
for any $r>0$, where $B$ is as above.
\\
As already observed in the Introduction, since
$\Phi \notin \Delta _{2}$ and $B\notin \Delta _{2}$, differentiability
of the functionals appearing in {\eqref{min_pbl}} is not guaranteed. Then
we cannot apply the standard method of Lagrange multipliers to obtain {Theorem~\ref{main_th}}. However, problem {\eqref{EL_bis}} can be still regarded as
the Euler-Lagrange equation associated with problem {\eqref{min_pbl}}.

Our next result guarantees the existence of a minimizer of problem {\eqref{min_pbl}}.

\begin{theorem}%
\label{teo3.3}
Under the same assumptions as in {Theorem~\ref{main_th}}, for any
$r>0$, minimization problem {\eqref{min_pbl}} has at least one minimizer
$u_{r}\in W^{1}_{0} L_{B, \Phi }(\Omega )$.
\end{theorem}

We observe that since no $\Delta _{2}$-condition is required on
$\Phi $ and on $B$, conditions
$\Phi _{\xi }\left ( \nabla u_{r}\right ) \in L_{\Phi _{\bullet }}(
\Omega ;\mathbb{R}^{N})$ and
$b(|u_{r}|) \in L_{B_{\bullet }}(\Omega )$ does not necessary occur, then
in general {\eqref{def_weaksol}} is not well-defined. Nevertheless, we are
able to prove the following result.

\begin{proposition}%
\label{prop2}
Under the same assumptions as in {Theorem~\ref{main_th}}, if $u_{r}\in W_{0}^{1}L_{B,
\Phi }(\Omega )$ is a minimizer of problem {\eqref{min_pbl}}, then
\begin{itemize}
\item[$(i)$]
$\Phi _{\xi }\left ( \nabla u_{r}\right ) \in L_{\Phi _{\bullet }}(
\Omega ;\mathbb{R}^{n})$;
\item[$(ii)$] $b(|u_{r}|) \in L_{B_{\bullet }}(\Omega )$.
\end{itemize}
\end{proposition}

\subsection{Examples}
\label{sec3.1}

In this Subsection, we specialize {Theorem~\ref{main_th}} to some classes of $N$-functions $\Phi $, which govern the differential
operator in the equation in {\eqref{EL_bis}}, with a distinctive structure.

\medskip

If $\Phi $ is defined as in {\eqref{ex_Ai}}, problem {\eqref{EL_bis}} takes
the form
%
\begin{equation}
\label{example22}
\begin{cases}
-\displaystyle \sum _{i=1}^{n} \left (A_{i}' (u_{x_{i}})\right )_{x_{i}}=
\lambda \, {b(|u|)} \, {\mathrm{sign }} \,u & \qquad \text{\textup{in\ }}
\Omega
\\
u=0 & \qquad \text{on}\;\partial \Omega \,,
\end{cases}
\end{equation}
where $A_{i}$, for $i=1,\ldots ,n$, are $1$-dimensional $N$-functions.
One has that (see \cite{C1})
%
\begin{equation}
\label{103}
\Phi _{\circ }(t) \approx \overline{A} (t) \qquad \text{near infinity},
\end{equation}
where $\overline{A}$ is the $1$-dimensional $N$-function obeying
%
\begin{equation}
\label{104}
\overline{A}^{-1} (\tau ) = \left ( \prod _{i=1}^{n} A_{i}^{-1} (
\tau )\right )^{\frac{1}{n}}\,.
\end{equation}
Thus, {Theorem~\ref{main_th}} holds for problem {\eqref{example22}} where $\Phi _{\circ }$ is replaced by $\overline{A}$ in
the definition of $\Phi _{n}$ (see {\eqref{sobolev_conj}}).
\bigskip

\noindent
\textbf{Example 1.} Let
%
\begin{equation}
\label{Ai}
A_{i}(t)= \frac{1}{p_{i}} t^{p_{i}} \log ^{\alpha _{i}}(c+ t) \quad
\text{ for } t>0,
\end{equation}
where $p_{i}>1$, $\alpha _{i} \in \Rvtex $, $i=1, \dots , n$, and $c$ positive
constant sufficiently large for all functions $A_{i}(t)$ be convex. Let
$\overline{p}$ and $\overline{\alpha }$ be defined as
\begin{equation*}
\label{p_bar}
\frac{1}{\overline{p}}= \frac{1}{n} \sum _{i=1}^{n} \frac{1}{p_{i}}\,
\quad \text{ and } \quad \overline{\alpha }= \frac{\overline{p}}{n}
\sum _{i=1}^{n} \frac{\alpha _{i}}{p_{i}}\,.
\end{equation*}
With this choice of $A_{i}$ in problem {\eqref{example22}} one has that function
$\Phi _{n} (x)$ defined as in {\eqref{sobolev_conj}} has the following behavior
\begin{equation*}
\label{202}
\Phi _{n} (x) \approx
\begin{cases}
t^{\overline{p}^{\star }} \left (\log (c + t)\right )^{
\frac{\overline{\alpha }n}{n-\overline{p}}}&\qquad \text{if}\;\;
\overline{p}<n
\\
e^{\frac{n}{n-1-\overline{\alpha }}}&\qquad \text{if}\;\; \overline{p}=n,
\, \overline{\alpha }< n-1
\\
e^{e^{\frac{n}{n-1}}} &\qquad \text{if}\;\; \overline{p}=n, \,
\overline{\alpha }= n-1\,,
\end{cases}
\end{equation*}
near infinity. When $\overline{p}>n$, or $\overline{p}=n$ and
$\overline{\alpha } >n-1$, condition {\eqref{int_finito}} holds.
\\
Assume that
\begin{align*}
\label{107}
B(t) \prec \prec
\begin{cases}
t^{\overline{p}^{\star }} \left (\log (c + t)\right )^{
\frac{\overline{\alpha }n}{n-\overline{p}}}&\qquad \text{if}\;\;
\overline{p}<n
\\
e^{\frac{n}{n-1-\overline{\alpha }}}&\qquad \text{if}\;\; \overline{p}=n,
\, \overline{\alpha }< n-1
\\
e^{e^{\frac{n}{n-1}}} &\qquad \text{if}\;\; \overline{p}=n, \,
\overline{\alpha }= n-1\,,
\end{cases}
\end{align*}
and
$B(t)\;\; \text{is any if }\;\;\overline{p} > n\;\; \text{ or}\; \;
\overline{p}=n, \,\overline{\alpha }> n-1$.
\\
Hence, thanks to {Theorem~\ref{main_th}}, for any $r>0$ there exist a constant
$\lambda _{r}>0$ and
$u_{r}\in W_{0}^{1}L^{\vec{p}}\log ^{\vec{\alpha }}L(\Omega )\cap L^{
\infty }(\Omega )$ such that $\int _{\Omega }B(u_{r})\, dx=r$ and
$u_{r}$ solves problem {\eqref{example22}} with of $A_{i}$ as in {\eqref{Ai}}, and with $\lambda =\lambda _{r}$. Here, $\vec{\alpha }$ stands
for the vector $(\alpha _{1}, \dots , \alpha _{n})$.
\bigskip

\noindent
\textbf{Example 2.} Now we show a special instance with
\begin{equation*}
\Phi (\xi )= \sum _{i=1}^{n} \left ( e^{|\xi _{i}|^{\alpha _{i}}} -1
\right )\qquad \text{for}\;\; \xi \in \rnvtex \,,
\end{equation*}
where $\alpha _{i} >1$.
\\
The corresponding problems read
\begin{equation*}
\label{example4}
\begin{cases}
-\displaystyle \sum _{i=1}^{n}\left ({\alpha _{i}}e^{|u_{x_{i}}|^{
\alpha _{i}}} |u_{x_{i}}|^{\alpha _{i}-2} u_{x_{i}} \right )_{x_{i}}=
\lambda \,{b(|u|)}\,{\mathrm{sign}}\, u & \qquad \text{\textup{in\ }} \Omega
\\
u=0 & \qquad \text{on}\;\partial \Omega \,.
\end{cases}
\end{equation*}
By  {\eqref{103}} and {\eqref{104}} again, we have
\begin{align*}
\label{118}
\Phi _{\circ }^{-1} (s)\approx \left (\prod _{i=1}^{n}\left ( \log (1+s)
\right )^{\frac{1}{\alpha _{i}}}\right )^{\frac{1}{n}}= \left (\log (1+s)
\right )^{\frac{1}{n} \sum _{i=1}^{n} \frac{1}{\alpha _{i}}} = \left (
\log (1+s)\right )^{\frac{1}{\overline{\alpha }}}\quad \text{near infinity}
\,,
\end{align*}
where $\overline{\alpha }$ is the harmonic average of $\alpha _{i}$, for
$i=1, \ldots , n$. Then,
\begin{equation*}
\label{119}
\Phi _{\circ }(t)\approx e^{t^{\overline{\alpha }}} -1 \qquad \text{near infinity}
\,,
\end{equation*}
and condition {\eqref{int_finito}} is always verified. Thus, {Theorem~\ref{main_th}} holds for any $N$-function $B$.
\bigskip

\noindent
\textbf{Example 3.} Let us consider now another particular case of the function {\eqref{ex_Ai}} given by
\begin{equation*}
\Phi (\xi )= \sum _{i=1}^{n-1} \frac{1}{p_{i}} |\xi _{i}|^{p_{i}} +
\left ( e^{|\xi _{n}|^{\alpha }} -1\right )\qquad \text{for
$\xi \in \rnvtex $}\,,
\end{equation*}
where $p_{i}>1$, for $i=1, \dots , n-1$, and $\alpha >1$. Note the
$\Phi \notin \Delta _{2}$. Now, problem {\eqref{example22}} agrees with
\begin{equation*}
\label{example3}
\begin{cases}
-\displaystyle \sum _{i=1}^{n-1}( |u_{x_{i}}|^{p_{i}-2} u_{x_{i}})_{x_{i}}
+ \left ( {\alpha } e^{|u_{x_{n}}|^{\alpha }} |u_{x_{n}}|^{\alpha -2}u_{x_{n}}
\right )_{x_{n}}= \lambda \,{b(|u|)}\,{\mathrm{sign}}\, u & \qquad \text{\textup{in\ }}
\Omega
\\
u=0 & \qquad \text{on}\;\partial \Omega \,.
\end{cases}
\end{equation*}
One can verify via {\eqref{103}} and {\eqref{104}} that
\begin{align*}
\label{108}
\Phi _{\circ }^{-1} (s) \approx \left (\prod _{i=1}^{n} A_{i}^{-1}(s)
\right )^{\frac{1}{n}}\approx \left (s^{\sum _{i=1}^{n-1} \frac{1}{p_{i}} }
\left (\log (1+s)\right )^{\frac{1}{\alpha }}\right )^{\frac{1}{n}}
\\
=s^{\frac{1}{n} \sum _{i=1}^{n-1} \frac{1}{p_{i}}} \left (\log (1+s)
\right )^{\frac{1}{n\alpha }} = s^{\frac{1}{\beta }} \left (\log (1+s)
\right )^{\frac{1}{n\alpha }}\quad \text{near infinity} \,,
\nonumber
\end{align*}
where $ \frac{1}{\beta } =\frac{1}{n} \sum _{i=1}^{n-1} \frac{1}{p_{i}}$. Then
\begin{equation*}
\label{109}
\Phi _{\circ }(t) \approx t^{\beta } \left (\log (1+t)\right )^{-
\frac{\beta }{n \alpha }}\quad \text{near infinity}\,.
\end{equation*}
If $\Phi _{\circ }$ verifies condition {\eqref{int_infinito}}, namely if
$ \frac{1-\beta }{n-1}> -1$, \emph{i.e.} if
$\sum _{i=1}^{n-1} \frac{1}{p_{i}}> 1$, then
\begin{equation*}
\label{114}
\Phi _{n}(t)\approx t^{\beta ^{*}} \left ( \log (1+t)\right )^{-
\frac{\beta }{\alpha (n-\beta )}}\qquad \text{near infinity}\,.
\end{equation*}
Whereas, if $\Phi $ fulfils condition {\eqref{int_finito}}, \emph{i.e.} if
$\sum _{i=1}^{n-1} \frac{1}{p_{i}}\leq 1$, then $\Phi _{n}$ agrees with
$+\infty $ near infinity.
\\
By assuming that
\begin{equation*}
\label{116}
B(t)\prec \prec s^{\beta ^{*}} \left ( \log (1+t)\right )^{-
\frac{\beta }{\alpha (n-\beta )}} \qquad \text{if}\;\; \sum _{i=1}^{n-1}
\frac{1}{p_{i}}> 1
\end{equation*}
and
\begin{equation*}
\label{117}
B(t) \qquad \text{is any} \qquad \text{if}\;\; \sum _{i=1}^{n-1}
\frac{1}{p_{i}}\leq 1\,,
\end{equation*}
{Theorem~\ref{main_th}} holds.
\bigskip

\noindent
\textbf{Example 4.} We present now a possible instance of examples which
generalize one from [Tr] provided by $N$-functions $\Phi $ of the form
\begin{equation*}
\label{ex_Tr}
\Phi \left ( \xi \right ) = \sum _{k=1}^{K} A_{k} \bigg (\Big |\sum _{i=1}^{n}
\alpha _{i}^{k} \xi _{i}\Big |\bigg )\qquad \text{for }\xi \in \rnvtex \,,
\end{equation*}
where $A_{k}$ are $N$-functions of one variable, $K\in \Nvtex $ and coefficients
$\alpha _{i}^{k}\in \Rvtex $ are arbitrary.
\\
When $n = 2$, we consider, for example, the $N$-function given by (see
\cite[Example 5]{ACCZ-G})
\begin{equation*}
\Phi (\xi ) = |\xi _{1} -\xi _{2}|^{p} + |\xi _{1}|^{q} (\log (c+ |
\xi _{1}|))^{\alpha }\qquad \text{for $\xi \in \mathbb{R}^{2}$}\,,
\end{equation*}
where $c$ is a sufficiently large constant for $\Phi $ to be convex,
$p>1$ and either $q\geq 1$ and $\alpha >0$, or $q=1$ and $\alpha >0$.
\\
Hence, problem {\eqref{EL_bis}} becomes
\begin{equation*}
\label{example5}
\begin{cases}
-\left [ \left (\Phi _{x_{1}} (u_{x_{1}}, u_{x_{2}})\right )_{x_{1}} +
\left (\Phi _{x_{2}} (u_{x_{1}}, u_{x_{2}})\right )_{x_{2}}\right ] =
\lambda \,{b(|u|)}\,{\mathrm{sign}}\, u & \quad \text{\textup{in\ }} \Omega
\\
u=0 & \quad \text{on}\;\partial \Omega \,,
\end{cases}
\end{equation*}
where
\begin{align*}
\Phi _{x_{1}} (u_{x_{1}}, u_{x_{2}})
&= p|u_{x_{1}} - u_{x_{2}}|^{p-2} (u_{x_{1}}
- u_{x_{2}}) + q |u_{x_{1}}|^{q-2} u_{x_{1}} \left (\log (c+ |u_{x_{1}}|)
\right )^{\alpha }
\\
&\quad {}+ \alpha |u_{x_{1}}|^{q-1}
\frac{u_{x_{1}}}{c+ |u_{x_{1}}|} \left (\log (c+ |u_{x_{1}}|)\right )^{
\alpha -1}\,,
\end{align*}
and
\begin{equation*}
\Phi _{x_{2}} (u_{x_{1}}, u_{x_{2}}) = - p |u_{x_{1}} - u_{x_{2}}|^{p-2}
(u_{x_{1}} - u_{x_{2}})\,.
\end{equation*}
The Sobolev conjugate of $\Phi $ takes the following values
%
\begin{equation}
\label{Sob_conj_Tr}
\Phi _{2} (s) \approx
\begin{cases}
s^{\frac{2pq}{p+q-pq}}\log ^{\frac{p\alpha }{p+q-pq}}(t)&\qquad \text{if}
\;\; pq<p+q
\\
\exp \big (t^{\frac{2(p+q)}{p+q-p\alpha }}\big )&\qquad \text{if}\;\; pq=p+q
\;\; \text{and}\;\; p\alpha <p+q
\\
\exp (\exp (t^{2})) &\qquad \text{if}\;\; pq=p\alpha =p+q
\end{cases}
\end{equation}
near infinity.
\\
If $pq=p+q$ and $p\alpha >p+q$, or $pq>p+q$, then condition {\eqref{int_finito}} holds.
\\
By assuming that
\begin{equation*}
\label{500}
B(t) \prec \prec \Phi _{2}(t)\,,
\end{equation*}
where $\Phi _{2}$ is as in {\eqref{Sob_conj_Tr}} and
\begin{equation*}
\label{501}
B(t) \; \; \text{is any} \qquad \text{if}\;\; pq=p+q\;\; \text{and}\;\; p
\alpha >p+q,\;\;\text{or}\;\; pq>p+q\,,
\end{equation*}
{Theorem~\ref{main_th}} holds.

\section{Proofs of main results}
\label{sec4}

In this section, we provide the proof of our main results stated in section
\S~3. In order to prove {Theorem~\ref{main_th}}, we first show the existence of
a minimizer of constrained minimization problem {\eqref{min_pbl}}. We only focus on the case $B\prec \prec \Phi _{n}$, since the other case runs
easily. In both cases, the key tool is the compact embedding of $W_{0}^{1}
L_{B,\Phi }(\Omega )$ into $E_{B}$ which is guaranteed by {Proposition~\ref{Prop_compact_embed}}.

\medskip
\textbf{Proof of {Theorem~\ref{teo3.3}.}}
\rm
Let us introduce the following functionals
$F: W_{0}^{1} L_{B,\Phi }(\Omega ) \vavtex \overline{\Rvtex }$ and
$G: W_{0}^{1} L_{B,\Phi }(\Omega ) \vavtex \overline{\Rvtex }$ defined as
\begin{equation*}
\label{F}
F(u)=\int _{\Omega }\Phi (\nabla u)\;\;dx\; ,
\end{equation*}
and
\begin{equation*}
\label{G}
G(u)=\int _{\Omega }B(u)\;\;dx\;,
\end{equation*}
respectively. We observe that $F$ is a finite-valued functional on
$W_{0}^{1} L_{B,\Phi }(\Omega )$ if and only if $\Phi $ fulfils the
$\Delta _{2}$-condition. Whereas, $G(u)$ is always finite for every
$u\in W_{0}^{1} L_{B,\Phi }(\Omega )$ because the compact embedding stated
in Proposition~\ref{Prop_compact_embed}.
\\
In order to prove the existence of a minimizer, we have to show first the
continuity of $G$ and the lower semicontinuity of $F$ with respect to the topology
$\sigma (W_{0}^{1} L_{B,\Phi }(\Omega ), W^{-1}E_{B_{\bullet },\Phi _{
\bullet }}(\Omega ))$, where $W_{0}^{1} L_{B,\Phi }(\Omega )$ has to be
understood as the dual space of
$W^{-1}E_{B_{\bullet },\Phi _{\bullet }}(\Omega )$.
\medskip

\noindent
\textbf{Step 1.}
\textit{$G$ is
$\sigma (W_{0}^{1} L_{B,\Phi }(\Omega ), W^{-1}E_{B_{\bullet },\Phi _{
\bullet }}(\Omega ))$ continuous.}

\noindent
It is enough to prove that if
%
\begin{equation}
\label{CD}
u_{h}\vavtex u\quad \text{in}\;\;\sigma (W_{0}^{1}L_{B,\Phi }(\Omega ), W^{-1}E_{B_{
\bullet },\Phi _{\bullet }}(\Omega )),
\end{equation}
then $G(u_{h}) \vavtex G(u)$.
\\
By {\eqref{CD}}, it follows that $u_{h}$ is bounded in
$W^{1}_{0}L_{B,\Phi }(\Omega )$. By the compact embedding of
$W_{0}^{1} L_{B,\Phi } (\Omega )$ in $E_{B}(\Omega )$ (see {Proposition~\ref{Prop_compact_embed}}), we have that $u_{h}$ converges to $u$ in norm
in $E_{B}(\Omega )$. Since convergence in norm implies the mean convergence
(see $\S $ 2.2), we get $B(2(u_{h}-u)) \rightarrow 0$ in
$L^{1}(\Omega )$. It follows, (up a subsequence), that
$u_{h} \rightarrow u$ \emph{a.e.} in $\Omega $ and there exists (up a subsequence)
a function $w \in L^{1}(\Omega )$ such that
\begin{equation*}
B(2|u_{h}-u|)\leq w(x) \qquad \text{ \emph{a.e.} in } \Omega \,.
\end{equation*}
Owing to the strictly monotonicity and the convexity of function $B$, we obtain
\begin{equation*}
B(u_{h}) \leq \frac{1}{2} B(2u) + \frac{1}{2}w \qquad
\text{ \emph{a.e.} in } \Omega \,,
\end{equation*}
and then the statement of Step 1 follows thanks to Lebesgue's dominate
convergence Theorem.

\medskip
\noindent
\textbf{Step 2.}
\textit{$F$ is
$\sigma (W_{0}^{1} L_{B,\Phi }(\Omega ), W^{-1}E_{B_{\bullet },\Phi _{
\bullet }}(\Omega ))$ lower semicontinuous.}

\noindent
By definition, it is enough to prove that $F(u)\leq \liminf F(u_{h})$ if {\eqref{CD}} holds.
\\
Let us fix $\varepsilon >0$. Since
$\nabla u_{h} \in L_{\Phi }(\Omega ; \rnvtex )$ for all $h$, by {Lemma~\ref{lemma_sup}}, there exists a function
$W \in E_{\Phi _{\bullet }}(\Omega ; \rnvtex )$ such that
\begin{equation*}
F(u_{h})=\int _{\Omega }\Phi (\nabla u_{h})\; dx \geq \int _{\Omega
}\nabla u_{h} \cdot W \; dx - \int _{\Omega }\Phi _{\bullet } (W)\; dx
\qquad \forall \; h\in \Nvtex
\end{equation*}
and
\begin{equation*}
F(u)=\int _{\Omega }\Phi (\nabla u)\; dx \leq \int _{\Omega }\nabla u
\cdot W \; dx - \int _{\Omega }\Phi _{\bullet } (W)\; dx + \varepsilon ,
\end{equation*}
namely
%
\begin{equation}
\label{Fn}
F(u_{h}) - F(u) \geq \int _{\Omega }\nabla u_{h} \cdot W \; dx - \int _{\Omega }\nabla u \cdot W \; dx - \varepsilon \qquad \forall \; h\in
\Nvtex .
\end{equation}
By {\eqref{Fn}} and {\eqref{CD}}, we get
\begin{equation*}
\liminf _{h} F(u_{h}) \geq F(u) - \varepsilon ,
\end{equation*}
and the proof of Step 2 follows by the arbitrariness of
$\varepsilon $.

\medskip
\noindent
\textbf{Step 3.} \textit{Existence of a minimizer of {\eqref{min_pbl}}.}

\noindent
Let $\{u_{h}\}_{h}\subset W^{1}_{0} L_{B,\Phi }(\Omega )$ be a minimizing
sequence of {\eqref{min_pbl}}, \emph{i.e.}
\begin{equation*}
G(u_{h})=\int _{\Omega }B(u_{h})\;\;dx =r \qquad \forall \; h\in \Nvtex
\end{equation*}
and
\begin{equation*}
F(u_{h})=\int _{\Omega }\Phi (\nabla u_{h})\;\; dx \rightarrow c_{r}
\qquad \text{as } h\rightarrow \infty .
\end{equation*}
This means that $\{\nabla u_{h}\}_{h}$ is bounded in mean and then in norm
in $L_{\Phi }(\Omega ;\rnvtex )$ (see $\S $ 2.2). By Banach-Alaoglu's Theorem,
there exists (up a subsequence)
$u_{r} \in W^{1}_{0} L_{B,\Phi }(\Omega ) $ such that
$ u_{h}\vavtex u_{r}$ in
$\sigma (W_{0}^{1}L_{B,\Phi }(\Omega ), W^{-1}E_{B_{\bullet },\Phi _{\bullet }}(\Omega ))$. By Step 1 and Step 2, it follows
\begin{equation*}
G(u_{r})=r \qquad \text{ and } \qquad F(u_{r})\leq \liminf F(u_{h})=c_{r}
\,.
\end{equation*}
By definition of $c_{r}$, we conclude that $F(u_{r})=c_{r}$.

\qed

Our next aim is to prove {Proposition~\ref{prop2}}. To do this, the next
auxiliary lemmas will be critical.

\begin{lemma}%
\label{Lemma_4.1}
Let $U\in \mathcal{L}_{\Phi }(\Omega ; \rnvtex )$. Then the following statements
hold
\begin{itemize}
\item[$(a)$]
$\Phi _{\xi }\left ( \left ( 1-\varepsilon \right ) U\right ) \in L_{
\Phi _{\bullet }}%
(\Omega ; \rnvtex )$ for all $\varepsilon \in \left ( 0,1\right ]$;
\item[$(b)$]
$\Phi _{\xi }\left ( \left ( 1-\varepsilon \right ) U+V\right ) \in L_{
\Phi _{\bullet }}(\Omega ; \rnvtex )$ for all
$V\in E_{\Phi }(\Omega ; \rnvtex )$ and for all
$\varepsilon \in \left ( 0,1\right ]$.
\end{itemize}
\end{lemma}

\medskip
\noindent \emph{Proof.}
\rm
Let $U\in \mathcal{L}_{\Phi }(\Omega ; \rnvtex )$ and
$V\in E_{\Phi }%
(\Omega ; \rnvtex )$. The case $\varepsilon =1$ is trivial, so let
$\varepsilon \in \left ( 0,1\right )$.
\\
\emph{$(a)$} By {(\ref{dis_gradiente})} and the convexity of $\Phi $, we get
\begin{align*}
\frac{\varepsilon }{1-\varepsilon }\,\Phi _{\bullet }\left ( \Phi _{
\xi }\left ( \left ( 1-\varepsilon \right ) U\right ) \right ) &\leq
\varepsilon \, U\cdot \Phi _{\xi }\left ( \left ( 1-\varepsilon
\right ) U\right )
\\
&\leq \varepsilon \, U\cdot \Phi _{\xi }\left ( \left ( 1-
\varepsilon \right ) U\right ) +\Phi \left ( \left ( 1-\varepsilon
\right ) U\right ) \leq \Phi \left ( U\right )\,.
\end{align*}
Then, since $U\in \mathcal{L}_{\Phi }(\Omega ; \rnvtex )$, it follows
$\Phi _{\bullet }\left ( \Phi _{\xi }\left ( \left ( 1-\varepsilon
\right ) U\right ) \right ) \in L^{1}(\Omega )$, namely
$\Phi _{\xi }\left ( \left ( 1-\varepsilon \right ) U\right ) \in
\mathcal{L}_{\Phi _{\bullet }}(\Omega ; \rnvtex )$.
\\
\emph{$(b)$} Thanks to the convexity of $\Phi $, we have that
%
\begin{align}
\label{304}
\Phi \left ( \frac{1}{1-\varepsilon /2}\left ( \left ( 1-\varepsilon
\right ) U+V\right ) \right ) \leq
\frac{1-\varepsilon }{1-\varepsilon /2}\Phi \left ( U\right ) +\left (
1-\frac{1}{1-\varepsilon /2}\right ) \Phi \left (
\frac{2V}{\varepsilon }\right )\,.
\end{align}
Inequality {\eqref{304}} gives
$\Phi \left ( \frac{1}{1-\varepsilon /2}\left ( \left ( 1-
\varepsilon \right ) U+V\right ) \right ) \in L^{1}(\Omega ,
\mathbb{R}^{n})$. Owing to \emph{$(a)$}, the statement \emph{$(b)$} follows.

\qed

For convenience of the reader, we state Lemma 4.2 in \cite{MT}.

\begin{lemma}
\label{Lemma_4.2}
Let $u,v\in E_{B}(\Omega )$, $u\nequiv 0$ and
$\int _{\Omega }B(u)v \; dx\neq 0$. Then the condition
\begin{equation*}
\int _{\Omega }B\left (\left ( 1-\varepsilon \right ) u+\delta v
\right )\; dx=\int _{\Omega }B(u)\; dx
\end{equation*}
defines a continuously differentiable function
$\delta =\delta (\varepsilon )$ in some interval
$\left ( -\varepsilon _{0},\varepsilon _{0}\right ) $ with
$\varepsilon _{0}>0$.
\\
Moreover, $\delta (0)=0$ and%
%
\begin{equation}
\label{delta_primo}
\delta ^{\prime }(0)=
\frac{\int _{\Omega }b(u)u\;
dx}{\int _{\Omega }b(u)v\; dx},
\end{equation}
where $b$ is the derivative of $B$.
\end{lemma}

\medskip
\textbf{Proof of {Proposition~\ref{prop2}.}}
\emph{$(i)$}. Let $r>0$
and $u_{r}\in W_{0}^{1}L_{B,\Phi }(\Omega )$ be a minimizer of {\eqref{min_pbl}}. Suppose, by contradiction, that $\Phi _{\xi }\left ( \nabla
u_{r}\right ) \notin L_{\Phi _{\bullet }}( \Omega ,\mathbb{R}^{N})$. By
{Proposition~\ref{equalityYoung}}, we get
%
\begin{equation}
\label{infinity}
\int _{\Omega }\Phi _{\xi }\left ( \nabla u_{r}\right ) \cdot \nabla u_{r}
\; dx =\int _{\Omega }\Phi (\nabla u_{r})\; dx+\int _{\Omega }\Phi _{
\bullet }(\Phi _{\xi }\left ( \nabla u_{r}\right ) )\; dx=+\infty .
\end{equation}
Now, we choose $v\in W_{0}^{1}E_{B,\Phi }(\Omega )$ such that
$\int _{\Omega } b(u_{r})u_{r}\, dx=\int _{\Omega }b(u_{r})v\, dx$. By
{Lemma~\ref{Lemma_4.2}} and by {\eqref{min_pbl}}, there exist
$\varepsilon _{0}\in (0, 1)$ and a function
$\delta \in \mathcal{C}^{1}(-\varepsilon _{0},\varepsilon _{0})$ fulfilling
%
\begin{equation}
\label{30}
\int _{\Omega }B(\left ( 1-\varepsilon \right ) u_{r}+\delta (
\varepsilon )v)\; dx=\int _{\Omega }B(u_{r})\; dx=r\qquad
\text{for all }\varepsilon \in (-\varepsilon _{0},\varepsilon _{0})\,.
\end{equation}
Moreover, $\delta (0)=0$ and $\delta '(0)=1$. Then, there exists
$\varepsilon _{1}\in (0, \varepsilon _{0})$ such that
$\delta (\varepsilon )\geq 0$ for all
$\varepsilon \in [0,\varepsilon _{1}]$ and
%
\begin{equation}
\label{200}
|\delta '(\varepsilon )|\leq \frac{3}{2} \qquad \text{for all }
\varepsilon \in [-\varepsilon _{1},\varepsilon _{1}]\,.
\end{equation}
By {\eqref{200}},
%
\begin{equation}
\label{201}
|\delta (\varepsilon )|\leq \frac{3}{2}\, |\varepsilon | \qquad
\text{for all } \varepsilon \in [-\varepsilon _{1},\varepsilon _{1}].
\end{equation}
Let us define the function $\Psi : [0, \varepsilon _{1}]\vavtex \Rvtex $ by
\begin{equation*}
\Psi (\varepsilon )=\int _{\Omega }\Phi (W_{\varepsilon }(x)) \; dx,
\end{equation*}
where
\begin{equation*}
\label{W}
W_{\varepsilon }(x)=\left ( 1-\varepsilon \right ) \nabla u_{r}(x)+
\delta (\varepsilon )\nabla v(x) \qquad \text{for}\; x\in \Omega .
\end{equation*}
Since, by {\eqref{dis_gradiente}},
$\Phi _{\xi }\left ( W_{\varepsilon }\right ) \cdot W_{\varepsilon }
\geq 0$ \textit{a.e.} in $\Omega $ and
$\Phi _{\xi }\left ( W_{\varepsilon }\right ) \cdot W_{\varepsilon }
\rightarrow \Phi _{\xi }\left ( \nabla u_{r}\right ) \cdot \nabla u_{r}$
\textit{a.e.} in $\Omega $, then by Fatou's Lemma and {\eqref{infinity}} we
have
%
\begin{equation}
\label{limite}
\int _{\Omega }\Phi _{\xi }\left ( W_{\varepsilon }\right ) \cdot W_{
\varepsilon }\; dx\rightarrow +\infty \qquad \text{ for } \varepsilon
\vavtex 0.
\end{equation}
Let $\varepsilon \in (0, 1)$. By the convexity of $\Phi $, estimates {\eqref{201}} and {\eqref{lin}}, it follows that
\begin{equation*}
\begin{split}%
\Phi (W_{\varepsilon })&\leq (1-\varepsilon ) \Phi (\nabla u_{r}) +
\varepsilon \Phi \Big (\frac{\delta (\varepsilon )}{\varepsilon }
\nabla v\Big )\leq \Phi (\nabla u_{r}) +
\frac{2\delta (\varepsilon )}{3} \Phi \Big (\frac{3}{2}\nabla v\Big )
\\
&\leq \Phi (\nabla u_{r}) + \Phi \Big (\frac{3}{2}\nabla v\Big )\in L^{1}(
\Omega )
\nonumber
\end{split}
\end{equation*}
for all $\varepsilon \in (0, \varepsilon _{1}]$, because
$v\in W_{0}^{1}E_{B,\Phi }(\Omega )$ and $u_{r}$ is a minimizer of problem {\eqref{min_pbl}}. Thanks to Lebesgue's dominated convergence Theorem,
$\Psi (\varepsilon )$ is continuous in $(0, \varepsilon _{1}]$. It is easily
to check the continuity in $\varepsilon =0$, as well.
\\
Let $\varepsilon _{2}\in (0, \varepsilon _{1})$ be arbitrary and set
$U= \displaystyle \frac{2}{2-\varepsilon }\left [(1-\varepsilon )
\nabla u_{r}+\delta (\varepsilon )\nabla v \right ]$. First, by {\eqref{dis_gradiente}} and convexity of $\Phi $ we get
%
\begin{align}
\label{15}
\Phi _{\bullet }\left (\Phi _{\xi }\left ((1-\varepsilon ) \nabla u_{r} +
\delta (\varepsilon ) \nabla v\right )\right )& = \Phi _{\bullet
}\left ( \Phi _{\xi }\left ( \left (1-\frac{\varepsilon }{2}\right )U
\right ) \right )\leq \Phi _{\xi }\left ( \left (1-
\frac{\varepsilon }{2}\right )U \right )\cdot \left (1-
\frac{\varepsilon }{2}\right )U
\\
&\leq \frac{2-\varepsilon }{\varepsilon }\left [ \Phi _{\xi }\left (
\left (1-\frac{\varepsilon }{2}\right )U\right ) \cdot
\frac{\varepsilon }{2}U+\Phi \left (\left (1-\frac{\varepsilon }{2}
\right )U\right ) \right ]
\nonumber
\\
&\leq \frac{2-\varepsilon }{\varepsilon }\Phi \left (U\right ).
\nonumber
\end{align}

Moreover, Young inequality {\eqref{Young}}, {\eqref{15}}, {\eqref{lin}}, {\eqref{200}}
and {\eqref{201}} yield
%
\begin{equation}
%
\begin{split}
\Big |\frac{\partial }{\partial \varepsilon }\Phi \left (W_{\varepsilon }\right )\Big | &=\frac{1}{\varepsilon }\Big |\Phi _{\xi
}\left (\left (1-\frac{\varepsilon }{2}\right )U\right ) \cdot \left (
\varepsilon \delta '(\varepsilon ) \nabla v-\varepsilon \nabla u_{r}
\right )\Big |
\\
\nonumber
&\leq \frac{1}{\varepsilon } \Big [\Phi _{\bullet }\left (\Phi _{\xi
}\left (\left (1-\frac{\varepsilon }{2}\right )U\right )\right )+
\Phi \left (\varepsilon \delta '(\varepsilon )\nabla v-\varepsilon
\nabla u_{r}\right ) \Big ]
\\
\nonumber
&\leq \frac{1}{\varepsilon } \Big [
\frac{2-\varepsilon }{\varepsilon }\Phi \left (
\frac{2-2\varepsilon }{2-\varepsilon } \nabla u_{r} +
\frac{2\delta (\varepsilon )}{2-\varepsilon } \nabla v \right ) +(1-
\varepsilon ) \Phi \left (
\frac{\varepsilon \delta '(\varepsilon )}{1-\varepsilon }\nabla v
\right )+\varepsilon \Phi \left ( \nabla u_{r}\right ) \Big ]
\\
\nonumber
&\leq \frac{1}{\varepsilon } \left [
\frac{2-\varepsilon }{\varepsilon }\left [ \left (
\frac{2-2\varepsilon }{2-\varepsilon }\right )\Phi \left ( \nabla u_{r}
\right ) + \left (\frac{\varepsilon }{2-\varepsilon }\right ) \Phi
\left (\frac{2\delta (\varepsilon )}{\varepsilon } \nabla v \right )
\right ]\right.
\\
&\quad {} \left.\vphantom{\left (\frac{2\delta (\varepsilon )}{\varepsilon } \nabla v \right )} + \varepsilon \delta '(\varepsilon ) \Phi \left (\nabla v
\right )+\varepsilon \Phi \left ( \nabla u_{r}\right ) \right ]
\\
\nonumber
&\leq \frac{1}{\varepsilon _{2}} \Big [ 2
\frac{1-\varepsilon _{2}}{\varepsilon _{2}}\Phi \left (\nabla u_{r}
\right ) + \frac{\delta (\varepsilon )}{2\varepsilon } \Phi \left (4
\nabla v \right ) +\frac{3}{2} \Phi \left (\nabla v\right )+\Phi
\left ( \nabla u_{r}\right ) \Big ]
\\
\nonumber
&\leq \frac{1}{\varepsilon _{2}} \Big [ 2
\frac{1-\varepsilon _{2}}{\varepsilon _{2}}\Phi \left (\nabla u_{r}
\right ) + \Phi \left (4\nabla v \right ) +\frac{3}{2} \Phi \left (
\nabla v\right )+\Phi \left ( \nabla u_{r}\right ) \Big ] \in L^{1}(
\Omega )
\end{split}
\end{equation}
for all $\varepsilon \in (\varepsilon _{2}, \varepsilon _{1})$, because
$v\in W_{0}^{1}E_{B,\Phi }(\Omega )$ and $u_{r}$ is a minimizer of problem {\eqref{min_pbl}}. Then we conclude that
%
\begin{equation}
\label{300}
\Psi '(\varepsilon )=\int _{\Omega }\Phi _{\xi }(W_{\varepsilon })\cdot (
\delta '(\varepsilon ) \nabla v - \nabla u_{r})\; dx
\end{equation}
for all $\varepsilon \in (\varepsilon _{2}, \varepsilon _{1})$. By the
arbitrariness of $\varepsilon _{2}$, equality {\eqref{300}} holds for all
$\varepsilon \in (0, \varepsilon _{1})$.
\\
We note that
%
\begin{align}
\label{110}
\Psi '(\varepsilon )&= \int _{\Omega }\Phi _{\xi }(W_{\varepsilon })
\left (\delta '(\varepsilon ) \nabla v -
\frac{W_{\varepsilon }- \delta ( \varepsilon ) \nabla v }{1-
\varepsilon }\right )\; dx
\\
&=- \frac{1}{1- \varepsilon } \int _{\Omega }\Phi _{\xi }(W_{\varepsilon })\cdot W_{\varepsilon }\; dx+ \left (\delta '(\varepsilon )+
\frac{\delta (\varepsilon )}{1-
\varepsilon }\right )\int _{\Omega }\Phi _{\xi }(W_{\varepsilon })\cdot
\nabla v\; dx \,.
\nonumber
\end{align}
Owing to Young inequality {\eqref{Young}} and inequality {\eqref{dis_gradiente}}, we have
%
\begin{equation}
\label{111}
\Phi _{\xi }(\xi ) \cdot 2\eta \leq \Phi _{\bullet }(\Phi _{\xi }(\xi ) )+
\Phi (2\eta ) \leq \xi \cdot \Phi _{\xi }(\xi ) + \Phi (2\eta )\quad
\forall \eta ,\xi \in \mathbb{R}^{n}.
\end{equation}
Since $v\in W_{0}^{1}E_{B,\Phi }(\Omega )$, by {\eqref{110}} and {\eqref{111}}, we can deduce that
%
\begin{align}
\label{4.2}
\Psi '(\varepsilon ) &\leq \left (\frac{\delta '
(\varepsilon )}{2}+\frac{\delta (\varepsilon ) /2-1}{1- \varepsilon }
\right )\int _{\Omega }\Phi _{\xi }(W_{\varepsilon })\cdot W_{\varepsilon
}\; dx+ \left (\frac{\delta '(\varepsilon )}{2}+
\frac{\delta (\varepsilon )}{2(1-
\varepsilon )}\right )\int _{\Omega }\Phi (2 \nabla v)\; dx
\\
&\leq \left (\frac{3}{4}+
\frac{\delta (\varepsilon ) /2-1}{1-
\varepsilon } \right )\int _{\Omega }\Phi _{\xi }(W_{\varepsilon })\cdot W_{\varepsilon }\; dx+ C \qquad \text{for}\;\;\varepsilon \in (0,
\varepsilon _{1}), \ \nonumber
\end{align}
where $C$ is a positive constant independent of $\varepsilon $. The last
estimate {(\ref{4.2})} and limit {\eqref{limite}} imply
$\lim _{\varepsilon \rightarrow 0} \Psi ^{\prime }(\varepsilon )=-
\infty $. Then, there exists $\varepsilon _{3}>0$ such that
$\Psi (\varepsilon _{3})<\Psi (0)$. On setting
$\widehat{u}_{r}=\left ( 1-\varepsilon _{3}\right ) u_{r}+\delta (
\varepsilon _{3})v$, we have
\begin{align*}
\int _{\Omega }\Phi \left ( \nabla \widehat{u}_{r}\right )\; dx =
\Psi (\varepsilon _{3})<\Psi (0)=\int _{\Omega }\Phi \left ( \nabla u_{r}
\right )\; dx
\end{align*}
and
\begin{equation*}
\int _{\Omega }B(\widehat{u}_{r})\; dx=r,
\end{equation*}
which is a contradiction recalling that $u_{r}$ is a minimizer. Hence,
$\Phi _{\xi }\left ( \nabla u_{r}\right ) \in \mathcal{L}_{\Phi _{
\bullet }}(\Omega ; \rnvtex )$ and the proof of \emph{(i)} is complete.

\emph{$(ii)$.} The idea of the proof is similar to that of
\emph{$(i)$}. For convenience of the reader, we give all details. Let
$r>0$ and let $u_{r}\in W_{0}^{1}L_{B,\Phi }(\Omega )$ be a minimizer of
problem {\eqref{min_pbl}}. Thanks to embedding {\eqref{compact_imbedding_1}}, $u_{r}\in E_{B}(\Omega )$. Let
$v\in E_{B}(\Omega )$ such that
$\int _{\Omega }%
b(u_{r})u_{r}\;dx=\int _{\Omega }b(u_{r})v\;dx$, and {Lemma~\ref{Lemma_4.2}} guaranties that there exist $\varepsilon _{0}\in (0,1)$ and a function
$\delta \in \mathcal{C}^{1}(-\varepsilon _{0},\varepsilon _{0})$ satisfying {\eqref{30}}. Moreover, $\delta (0)=0$, $ \delta (\varepsilon )\geq 0$ for
all $\varepsilon \in [0,\varepsilon _{1}]$, and {\eqref{200}} and {\eqref{201}} hold. On setting
\begin{equation*}
\label{35}
\Lambda (\varepsilon )=\int _{\Omega }B(\omega _{\varepsilon }(x)) \; dx
\qquad \text{with}\;\;\varepsilon \in [0, \varepsilon _{1}]\,,
\end{equation*}
where
$\omega _{\varepsilon }(x) = (1- \varepsilon ) u_{r}(x) + \delta (
\varepsilon ) v(x)$ for $x\in \Omega $, by {\eqref{30}}, it follows that
$ \Lambda (\varepsilon ) = r $ and then
$ \Lambda '(\varepsilon ) = 0$ for every
$\varepsilon \in [0, \varepsilon _{1}]$. Now, we assume by absurdum that
$b(u_{r})\notin L_{B_{\bullet }}(\Omega )$, \emph{i.e.}
%
\begin{equation}
\label{777}
\int _{\Omega }B_{\bullet }(b(u_{r}))\ dx = + \infty \,.
\end{equation}
Let $\varepsilon _{2}\in (0,\varepsilon _{1})$ be arbitrary. The monotonicity
of $b$ and {Lemma~\ref{Lemma_4.1}} for $1$-dimensional Young function give
\begin{equation*}
\label{38}
\begin{split} \Big |\frac{\partial }{\partial \varepsilon } B(\omega _{\varepsilon }(x))\Big | &= \Big | b\left ((1- \varepsilon ) u_{r}(x) +
\delta (\varepsilon ) v(x)\right )(\delta '(\varepsilon ) v(x)- u_{r}(x))
\Big |
\\
&\leq \Big | b\left ((1- \varepsilon ) u_{r}(x) + \delta (\varepsilon )
v(x)\right )\Big | \Big ( |\delta '(\varepsilon )| |v(x)| + |u_{r}(x)|
\Big )
\nonumber
\\
& \leq \Big | b\Big ((1- \varepsilon _{2})\, u_{r}(x) + \frac{3}{2} v(x)
\Big )\Big |\Big ( \frac{3}{2} |v(x)| + |u_{r}(x)|\Big )\in L^{1}(\Omega )
\nonumber
\end{split}
\end{equation*}
for any $ \varepsilon \in (\varepsilon _{2}, \varepsilon _{1})$. For the
arbitrariness of $\varepsilon _{2}$, it follows
\begin{equation*}
\Lambda '(\varepsilon )= \int _{\Omega }b(\omega _{\varepsilon })(\delta '(
\varepsilon ) v(x)- u_{r}(x))\; dx
\end{equation*}
for every $\varepsilon \in (0, \varepsilon _{1})$.
\\
By  estimate $b(s)t\leq \frac{1}{2} b(s) s + b(2t)t$ for all
$s,t \in \Rvtex $, we have
%
\begin{align}
\label{39}
\Lambda '(\varepsilon ) &= \int _{\Omega }b(\omega _{\varepsilon })\left (
\delta '(\varepsilon ) v(x) +
\frac{\delta (\varepsilon ) v(x) - \omega _{\varepsilon }}{1-\varepsilon }
\right )\; dx
\\
& = -\frac{1}{1-\varepsilon } \int _{\Omega }b(\omega _{\varepsilon })
\omega _{\varepsilon }\; dx+ \int _{\Omega }\left (\delta '(\varepsilon ) +
\frac{\delta (\varepsilon )}{1-\varepsilon }\right )b(\omega _{\varepsilon }) v(x) \; dx
\nonumber
\\
& \leq \left (\frac{\delta ' (\varepsilon ) }{2} +
\frac{1}{1-\varepsilon } \left ( \frac{\delta (\varepsilon )}{2} -1
\right )\right ) \int _{\Omega }b(\omega _{\varepsilon }) \omega _{\varepsilon }\; dx +\left (\delta '(\varepsilon ) +
\frac{\delta (\varepsilon )}{1-\varepsilon }\right )\int _{\Omega }b(2 v)
v \; dx\ \nonumber
\end{align}
for every $\varepsilon \in (0, \varepsilon _{1})$. Young inequality and {\eqref{dis_gradiente}} yield
%
\begin{equation}
\label{40}
\int _{\Omega }b(2 v) v\; dx \leq \int _{\Omega }B(v)\; dx +
\int _{\Omega }B_{\bullet }(b(2v))\; dx \leq \int _{\Omega }B(v)\; dx + \int _{\Omega }B ( 4v)\; dx < +\infty \,.
\end{equation}
By {Proposition~\ref{equalityYoung}} and {\eqref{777}}, we get
%
\begin{equation}
\label{41}
\int _{\Omega }b(u_{r}) u_{r} \; dx = \int _{\Omega }B(u_{r})\; dx +
\int _{\Omega }B_{\bullet }(b(u_{r}))\; dx = +\infty \,.
\end{equation}
By the continuity of $b$, it follows that
$ b(\omega _{\varepsilon }) \omega _{\varepsilon }\vavtex b(u_{r}) u_{r}$
\emph{a.e.} in $\Omega $. Since
$b(\omega _{\varepsilon }) \omega _{\varepsilon }\geq 0$ \emph{a.e.} in
$\Omega $, Fatou's Lemma and {\eqref{41}} yield
%
\begin{equation}
\label{43}
\lim _{\varepsilon \vavtex 0^{+}} \int _{\Omega }b(w_{\varepsilon }) w_{\varepsilon }= \int _{\Omega }b(u_{r}) u_{r}\; dx =+\infty \,.
\end{equation}
Now, we pass to the limit in the right-hand side of {\eqref{39}}. By {\eqref{40}} and {\eqref{43}}, we have
$ \lim _{\varepsilon \vavtex 0^{+}} \Lambda '(\varepsilon )= - \infty $ that
is in contradiction with the fact that $\Lambda '(\varepsilon )=0$.
\\
This implies that $b(u_{r})\in \mathcal{L}_{B_{\bullet }}(\Omega )$ and
complete the proof of \emph{$(ii)$}.

\qed

We are now in a position to accomplish the proof of {Theorem~\ref{main_th}}. In what follows it is important to underline that
$\big (W^{1}_{0} L_{B,\Phi }(\Omega ), W^{1}_{0} E_{B,\Phi }(\Omega ); W^{-1}
L_{B_{\bullet },\Phi _{\bullet }}(\Omega )$,
$W^{-1}E_{b_{\bullet },\Phi _{\bullet }}(\Omega )\big )$ is a complementary
system (see Subsection \ref{sec_compl_syst}) under our assumption on
$\Omega $.

\textbf{Proof of {Theorem~\ref{main_th}.}}
\rm
Let us define the functionals $dF$ ad $dG$ by
\begin{equation*}
\label{dF}
\langle dF, v\rangle = \int _{\Omega }\Phi _{\xi }(\nabla u_{r})\cdot
\nabla v \; dx
\end{equation*}
and
\begin{equation*}
\label{dG}
\langle dG, v\rangle = \int _{\Omega }\frac{b(|u_{r}|)}{|u_{r}|}u_{r} v
\; dx
\end{equation*}
for any $v\in W_{0}^{1} E_{B,\Phi }(\Omega )$, where $u_{r}$ is a minimizer
of problem {\eqref{min_pbl}}. By {Proposition~\ref{prop2}}, the previous functionals
are well-defined. Set
\begin{equation*}
{\mathrm{Ker}}\;dF=\{ v\in W_{0}^{1}E_{B,\Phi } (\Omega ): \langle dF, v
\rangle =0\}
\end{equation*}
and
\begin{equation*}
{\mathrm{Ker}}\;dG=\{ v\in W_{0}^{1}E_{B,\Phi } (\Omega ): \langle dG, v
\rangle =0\}\,.
\end{equation*}
If we prove that
%
\begin{equation}
\label{kern_inclusion}
{\mathrm{Ker}}\;dF\subset {\mathrm{Ker}}\;dG \,,
\end{equation}
then Proposition 43.1 in \cite{Ze} assures the existence of
$\lambda _{r}\in \Rvtex $, associated with the minimizer $u_{r}$, such that
equality {\eqref{def_sol}}, for $u=u_{r}$, holds for any test function
$\varphi $ in $W_{0}^{1}E_{B,\Phi }(\Omega )$. Finally, the
$\sigma (W_{0}^{1}L_{B,\Phi }(\Omega ), W^{-1}L_{B_{\bullet },\Phi _{\bullet }}(\Omega ))$-density of $W_{0}^{1}E_{B,\Phi }(\Omega )$ in
$W_{0}^{1}L_{B,\Phi }(\Omega )$ (see {Lemma~\ref{lemma_E_dense}} above) guarantees
that {\eqref{kern_inclusion}} it is enough to conclude.
\\
Then our goal is to prove {\eqref{kern_inclusion}}, which will follow by
the inclusion
\begin{align*}
V_{G} &:=\left \{  v\in W_{0}^{1}E_{B,\Phi }(\Omega ): \int _{\Omega
}\frac{b(|u_{r}|)}{|u_{r}|} u_{r} v \; dx >0\right \}
\\
&\subset \left \{  v
\in W_{0}^{1}E_{B,\Phi }(\Omega ): \int _{\Omega }\Phi _{\xi }(\nabla u_{r})
\cdot \nabla v \; dx >0\right \}  :=V_{F}\,.
\end{align*}
In order to verify the last inclusion, let us consider an arbitrary
$v\in V_{G}$. By {Lemma~\ref{Lemma_4.2}}, there exist
$\varepsilon _{0}\in (0, 1)$ and
$\delta \in \mathcal{C}^{1}(-\varepsilon _{0}, \varepsilon _{0})$ such
that
%
\begin{equation}
\label{26}
\int _{\Omega }B((1- \varepsilon ) u_{r} + \delta (\varepsilon ) v)\; dx
=\int _{\Omega }B(u_{r})\; dx=r \qquad \forall \varepsilon \in (-
\varepsilon _{0}, \varepsilon _{0}).
\end{equation}
On setting
$w_{\varepsilon }=(1-\varepsilon ) u_{r} + \delta (\varepsilon )v$, the definition
of $u_{r}$ and {\eqref{26}} assure that
%
\begin{equation}
\label{20}
\int _{\Omega
}\frac{\Phi (\nabla w_{\varepsilon }) - \Phi (\nabla u_{r})}{\delta (\varepsilon )}
\;dx \geq 0 \qquad \forall \varepsilon \in (0,\varepsilon _{1}).
\end{equation}
By {\eqref{delta_primo}}, we have $\delta '(0)>0$ and there exists
$\varepsilon _{1} \in (0, \varepsilon _{0})$ such that
$\frac{\delta '(0)}{2}<\delta '(\varepsilon )< 2\delta '(0)$ for all
$\varepsilon \in (-\varepsilon _{1}, \varepsilon _{1})$. By integrating
with respect to $\varepsilon $, we obtain
%
\begin{equation}
\label{301}
\frac{\delta '(0)}{2}<\frac{\delta (\varepsilon )}{\varepsilon }< 2
\delta '(0) \qquad \forall \varepsilon \in (0, \varepsilon _{1})\,.
\end{equation}
Since $\Phi $ is locally Lipschitz, we get in the ball $B(0,R) $ of radius
$R=|\nabla u_{r}|+\frac{3}{2}|\nabla v|$ that
%
\begin{align}
\label{23}
\bigg |
\frac{\Phi (\nabla w_{\varepsilon }) - \Phi (\nabla u_{r})}{\delta (\varepsilon )}
\bigg | &\leq \frac{L}{\delta (\varepsilon )} \left | \nabla w_{\varepsilon }- \nabla u_{r} \right | \leq L
\frac{\varepsilon }{\delta (\varepsilon )}|\nabla u_{r}|+L|\nabla v|
\leq L \frac{2}{\delta '(0)} |\nabla u_{r}| + L|\nabla v|,
\end{align}
where $L$ is the Lipschitz constant and the last inequality follow by {\eqref{301}}. The rightmost side of {\eqref{23}} belongs to
$L^{1}(\Omega )$ because $v \in W_{0}^{1}E_{B,\Phi }(\Omega )$ and
$u_{r}$ is the solution to {\eqref{min_pbl}}.

On recalling that $\Phi $ is differentiable, easily computation gives
\begin{equation*}
\lim _{\varepsilon \vavtex 0^{+}}
\frac{\Phi (\nabla w_{\varepsilon }) - \Phi (\nabla u_{r})}{\delta (\varepsilon )}=
\Phi _{\xi }(\nabla u_{r})\cdot \nabla v - \Phi _{\xi }(\nabla u_{r})
\cdot \frac{\nabla u_{r}}{\delta '(0)} \quad \text{ \textit{a.e.} in }
\Omega .
\end{equation*}
Then, by Lebesgue's dominate convergence Theorem, it follows
%
\begin{equation}
\label{24}
\lim _{\varepsilon \vavtex 0^{+}} \int _{\Omega
}\frac{\Phi (\nabla w_{\varepsilon }) - \Phi (\nabla u_{r})}{\delta (\varepsilon )}
\; dx = \int _{\Omega }\left [ \Phi _{\xi }(\nabla u_{r})\cdot \nabla v -
\Phi _{\xi }(\nabla u_{r})\cdot \frac{\nabla u_{r}}{\delta '(0)} \; dx
\right ]\, .
\end{equation}
Combining {\eqref{20}} and {\eqref{24}} we have
\begin{equation*}
\label{25}
\int _{\Omega }\Phi _{\xi }(\nabla u_{r}) \cdot \nabla v\; dx \geq
\frac{1}{\delta '(0)}\int _{\Omega }\Phi _{\xi }(\nabla u_{r})\cdot \nabla u_{r}
\; dx>0,
\end{equation*}
namely $v\in V_{F}$. Then $V_{G}\subset V_{F}$ follows by arbitrariness
of $v\in V_{F}$.

At this point of the proof, we have found a function $u_{r}$ in
$W_{0}^{1}L_{B, \Phi }(\Omega )$ such that $\int _{\Omega }B(u_{r})\; dx=r$,
and fulfils the equation
%
\begin{equation}
\label{sol_weak}
\int _{\Omega }\Phi _{\xi }(\nabla u_{r}) \cdot \nabla \varphi \; dx =
\lambda _{r} \int _{\Omega }\frac{b(|u_{r}|)}{|u_{r}|}u_{r} \varphi \; dx
\end{equation}
for any $\varphi \in W_{0}^{1}L_{B,\Phi }(\Omega )$.
\\
Our aim is now to prove that this function $u_{r}$ is actually a weak solution
to problem {\eqref{EL_bis}} as stated in {Definition~\ref{def_weaksol}}.
\\
To do this, we first observe that, by inclusion {\eqref{W_subset_W}},
$ u_{r}\in \mathcal{W}_{0}^{1}L_{\Phi }(\Omega )$.
\\
Next, by \cite[Proposition 2.4]{ACCZ-G}, one has that, given any function
$\varphi \in \mathcal{W}_{0}^{1} L_{\Phi }(\Omega )\cap L^{\infty }(
\Omega )$, there exist a constant $C=C(\Omega )$ and a sequence
$\{\varphi _{k}\}_{k} \subset \mathcal{C}_{0}^{\infty }(\Omega )$ such that
%
\begin{equation}
\label{sep30}
\varphi _{k} \tovtex \varphi \quad \text{a.e. in $\Omega $,}
\end{equation}
%
\begin{equation}
\label{july36}
\|\varphi _{k}\|_{L^{\infty }(\Omega )}\leq C\|\varphi \|_{L^{\infty }(
\Omega )} \quad \text{for every $k\in \Nvtex $,}
\end{equation}
%
\begin{equation}
\label{angela1grad}
\nabla \varphi _{k} \tovtex \nabla \varphi \quad \text{modularly in
$L_{\Phi }(\Omega ;\rnvtex )$.}
\end{equation}
Then, we have
%
\begin{equation}
\label{sol_weak_k}
\int _{\Omega }\Phi _{\xi }(\nabla u_{r}) \cdot \nabla \varphi _{k} \; dx =
\lambda _{r} \int _{\Omega }\frac{b(|u_{r}|)}{|u_{r}|}u_{r} \varphi _{k}
\; dx
\end{equation}
for any $\varphi _{k}\in \mathcal{C}_{0}^{\infty }(\Omega )$.
\\
Condition {\eqref{angela1grad}} means that there exists a constant
$l>0$ such that
\begin{equation*}
\lim _{k\vavtex \infty } \int _{\Omega } \Phi \left (
\frac{\nabla \varphi _{k} - \nabla \varphi }{l}\right )\; dx =0\,.
\end{equation*}
Moreover, \cite[Proposition 2.2]{ACCZ-G} yields that, owing to conditions {\eqref{angela1grad}} and
$\Phi _{\xi }(\nabla u_{r})\,{\in}\, L_{\Phi _{\bullet }}(\Omega {;} \rnvtex )$, there
exists a subsequence of $\{\nabla \varphi _{k}\}_{k}$, still indexed by
$k$, such that
\begin{equation*}
\lim _{k\tovtex \infty }\int _{\Omega }\Phi _{\xi }(\nabla u_{r}) \cdot
\nabla \varphi _{k} \; dx = \int _{\Omega }\Phi _{\xi }(\nabla u_{r})
\cdot \nabla \varphi \; dx\,.
\end{equation*}
Therefore, by the dominated convergence theorem coupling with {\eqref{sep30}}, {\eqref{july36}} and the fact that
$b(|u_{r}|)\in L_{B_{\bullet }}(\Omega )$, and hence in
$L^{1}(\Omega )$, we have that
\begin{equation*}
\lim _{k\tovtex \infty } \int _{\Omega }\frac{b(|u_{r}|)}{|u_{r}|}u_{r}
\varphi _{k} \; dx = \int _{\Omega }\frac{b(|u_{r}|)}{|u_{r}|}u_{r}
\varphi \; dx\,.
\end{equation*}
Then, $u_{r}$ is a weak solution to {\eqref{EL_bis}} in $\mathcal{W}_{0}^{1}L_{
\Phi }(\Omega )$ in the sense of {Definition~\ref{def_weaksol}}. Concerning the
boundedness of $u_{r}$, it follows directly by \cite[Theorem 4.1]{A}.

\qed

\paragraph*{Funding} This research was partly supported by: (i) GNAMPA (National
Group for Mathematical Analysis, Probability and their Applications) of
the Italian INdAM - National Institute of High Mathematics (Grant No. not available); (ii) the third
author is additionally partly funded by FFABR (Grant No. not available).


\end{document}